\documentclass[11pt]{article}
\usepackage[%
tmargin=1in,bmargin=1in,%
lmargin=1.4in,rmargin=1.4in,%
marginparsep=0.2in, marginparwidth=1.0in%
]{geometry}
\AtBeginDocument{%
	\setlength{\abovedisplayskip}{1.5ex plus 0.3ex minus 0.3ex}%
	\setlength{\abovedisplayshortskip}{1.0ex plus 0.3ex minus 0.3ex}%
	\setlength{\belowdisplayskip}{1.5ex plus 0.3ex minus 0.3ex}%
	\setlength{\belowdisplayshortskip}{1.0ex plus 0.3ex minus 0.3ex}%
}
\usepackage{amscd}
\usepackage{amsfonts}
\usepackage{amsmath}
\usepackage{amssymb}
\usepackage{amsthm}
\usepackage{fancyhdr}
\usepackage{latexsym}
\usepackage{xparse}
\usepackage{physics}
\usepackage{xspace}
\usepackage{titlesec}
\usepackage[dotinlabels]{titletoc}
\usepackage{textcase}
\usepackage[colorlinks=true, pdfstartview=FitV, linkcolor=black,
citecolor=black, urlcolor=blue]{hyperref}
\usepackage[style=alphabetic,maxbibnames=99,backend=bibtex8]{biblatex}
\usepackage{relsize}
\usepackage[bbgreekl]{mathbbol}

\usepackage{ulem}

\usepackage{cancel}

\DeclareSymbolFontAlphabet{\mathbb}{AMSb}
\DeclareSymbolFontAlphabet{\mathbbl}{bbold}

\usepackage{color}

\definecolor{Red}{rgb}{1.00, 0.00, 0.00}
\definecolor{DarkGreen}{rgb}{0.00, 1.00, 0.00}
\definecolor{Blue}{rgb}{0.00, 0.00, 1.00}
\definecolor{Cyan}{rgb}{0.00, 1.00, 1.00}
\definecolor{Magenta}{rgb}{1.00, 0.00, 1.00}
\definecolor{DeepSkyBlue}{rgb}{0.00, 0.75, 1.00}
\definecolor{DarkGreen}{rgb}{0.00, 0.39, 0.00}
\definecolor{SpringGreen}{rgb}{0.00, 1.00, 0.50}
\definecolor{DarkOrange}{rgb}{1.00, 0.55, 0.00}
\definecolor{OrangeRed}{rgb}{1.00, 0.27, 0.00}
\definecolor{DeepPink}{rgb}{1.00, 0.08, 0.57}
\definecolor{DarkViolet}{rgb}{0.58, 0.00, 0.82}
\definecolor{SaddleBrown}{rgb}{0.54, 0.27, 0.07}
\definecolor{Black}{rgb}{0.00, 0.00, 0.00}
\definecolor{dark-magenta}{rgb}{.5,0,.5}
\definecolor{myblack}{rgb}{0,0,0}
\definecolor{darkgray}{gray}{0.5}
\definecolor{lightgray}{gray}{0.75}

\usepackage[pdftex]{graphicx}
\usepackage{epstopdf}

\newcommand{\zz}{\mathbb{Z}}

\newcommand{\bb}[1]{\mathbb{#1}}

\newcommand{\mc}[1]{\mathcal{#1}}

\newcommand{\bs}{\backslash}
\newcommand{\fff}{\mathbb{F}}
\newcommand{\simto}{\overset{\sim}{\longrightarrow}}

\DeclareMathOperator{\colim}{colim}

\DeclareMathOperator{\Spec}{Spec}

\newcommand{\mo}[1]{\operatorname{#1}}

\newcommand{\f}[2]{\frac{#1}{#2}}

\newcommand{\bdot}{\text{\textbullet}}
\newcommand{\syn}{\mo{syn}}
\newcommand{\tors}{\mo{tors}}
\newcommand{\mf}[1]{\mathfrak{#1}}
\newcommand{\cat}[1]{\mo{Arith}^{S^1}(#1)}
\newcommand{\Prism}{{\mathlarger{\mathbbl{\Delta}}}}
\newcommand{\fge}{\mo{F}\text{-}\mo{Gauge}_{\Prism}^{\mo{perf}}}
\newcommand{\lge}{\mo{Perf}\left((\fff_p)_{\acute{e}t},\zz_{\ell}\right)}
\DeclareUnicodeCharacter{202F}{\,}
\newtheorem*{short hand}{theorem name}

\def\makeautorefname#1#2{\expandafter\def\csname#1autorefname\endcsname{#2}}
\def\equationautorefname~#1\null{(#1)\null}
\makeautorefname{footnote}{footnote}%
\makeautorefname{item}{item}%
\makeautorefname{figure}{Figure}%
\makeautorefname{table}{Table}%
\makeautorefname{part}{Part}%
\makeautorefname{appendix}{Appendix}%
\makeautorefname{chapter}{Chapter}%
\makeautorefname{section}{Section}%
\makeautorefname{subsection}{Section}%
\makeautorefname{subsubsection}{Section}%
\makeautorefname{theorem}{Theorem}%
\makeautorefname{corollary}{Corollary}%
\makeautorefname{lemma}{Lemma}%
\makeautorefname{proposition}{Proposition}%
\makeautorefname{pro}{Property}
\makeautorefname{conj}{Conjecture}%
\makeautorefname{definition}{Definition}%
\makeautorefname{notn}{Notation}
\makeautorefname{notns}{Notations}
\makeautorefname{rem}{Remark}%
\makeautorefname{quest}{Question}%
\makeautorefname{example}{Example}%
\theoremstyle{plain}  
\newtheorem{theorem}{Theorem}[section]
\newtheorem{proposition}{Proposition}[section]
\newtheorem{lemma}{Lemma}[section]
\newtheorem{corollary}{Corollary}[section]
\newtheorem{conjecture}[subsection]{Conjecture}

\theoremstyle{definition}
\newtheorem{definition}{Definition}[section]

\newtheorem{remark}{Remark}[section]

\makeatletter
\let\c@corollary=\c@theorem
\let\c@proposition=\c@theorem
\let\c@lemma=\c@theorem
\let\c@definition=\c@theorem
\let\c@example=\c@theorem
\let\c@remark=\c@theorem
\numberwithin{equation}{section}
\makeatother

\usepackage{tikz}

\usetikzlibrary{patterns}
\usepackage{tikz-cd}
\usepackage{spectralsequences}
\SseqOrientationSideways

\usepackage[T1]{fontenc}
\usepackage[largesc]{newtxtext}
\usepackage{newtxmath}

\usepackage[%
cal=cm,       calscaled=0.95,%
bb=boondox,   bbscaled=1.0,%
frak=euler,   frakscaled=1.05,%
]{mathalfa}  

\makeatletter
\g@addto@macro\bfseries{\boldmath}
\makeatother
\usepackage[textsize=footnotesize,backgroundcolor=orange!50]{todonotes}
\usepackage[compress,capitalize]{cleveref}

\let\theoldbibliography\thebibliography
\renewcommand{\thebibliography}[1]{%
	\theoldbibliography{#1}%
	\setlength{\parskip}{0ex}
	\setlength{\itemsep}{0.5ex plus 0.2ex minus 0.2ex}
	\small
}
\apptocmd{\thebibliography}{\raggedright}{}{}

\renewcommand{\title}[1]{\newcommand{\thetitle}{#1}}
\renewcommand{\author}[1]{\newcommand{\theauthor}{#1}}

\renewcommand{\maketitle}{%
	\begin{center}
		{\linespread{1.15}%
			\bfseries\MakeTextUppercase%
			\thetitle\par} \vspace{4.0ex}
		\footnotesize
		{\MakeUppercase \theauthor}
	\end{center}
	\vspace{3.0ex}
	\thispagestyle{fancy}
}

\renewenvironment{abstract}{\noindent\begin{center}\begin{minipage}{0.8\linewidth}\small{\scshape Abstract.}}{\end{minipage}\end{center}}
\newlength{\tagsep}
\makeatletter
\def\fullwidthdisplay{\displayindent\z@ \displaywidth\columnwidth}
\edef\@tempa{\noexpand\fullwidthdisplay\the\everydisplay}
\everydisplay\expandafter{\@tempa}
\makeatother

\titleformat{\section}{\centering}{\textsection\thesection.}{1.5\tagsep}{\scshape}
\titleformat{\subsection}[runin]{}{\fontseries{b}\selectfont\textsection\bfseries\thesubsection.}{1.5\tagsep}{\bfseries}[.]

\titlespacing*{\section}{0pt}{4ex}{\medskipamount}
\titlespacing*{\subsection}{0pt}{\bigskipamount}{0.5em}

\dottedcontents{section}[1.3em]{\vspace{0.25ex}}{1.3em}{0.7em}
\dottedcontents{subsection}[2.0em]{\vspace{0.25ex}}{2.0em}{0.7em}

\title{
Special Values without Semi-simplicity Via K-Theory
}

\author{ Logan Hyslop }

\date{mm/dd/yyyy}

\bibliography{Bibliography.bib}

\begin{document}
\DeclareDocumentCommand\rart{ g }{%
	{\ar[r, tail]%
		\IfNoValueF {#1} { \ar[r, tail, "#1"]}%
	}%
}
\DeclareDocumentCommand\dart{ g }{%
	{\ar[d, tail]%
		\IfNoValueF {#1} { \ar[d, tail, "#1"]}%
	}%
}
\DeclareDocumentCommand\rarh{ g }{%
	{\ar[r, two heads]%
		\IfNoValueF {#1} { \ar[r, two heads, "#1"]}%
	}%
}
\DeclareDocumentCommand\darh{ g }{%
	{\ar[d, two heads]%
		\IfNoValueF {#1} { \ar[d, two heads, "#1"]}%
	}%
}

\maketitle



%
%
%
%
%
%
%
\begin{abstract}
In this paper, motivated by studying special values of zeta functions attached to finite type $\fff_p$-schemes, we introduce a category $\cat{R}$ of ``arithmetic $C(S^1,R)$-modules'' attached to any Dedekind ring $R$, and compute the 0th $K$-group $K_0(\cat{R})$ of this category.  Specializing to the case of $R=\zz_\ell$ for some prime $\ell\neq p$ (resp. $R=\zz_p$), we prove that there is a natural functorial lift of the \'{e}tale cohomology of perfect \'{e}tale $\zz_{\ell}$ sheaves (resp. syntomic cohomology of perfect prismatic F-gauges) on a point to $\cat{\zz_{\ell}}$ (resp. $\cat{\zz_p}$).  This allows us to define a notion of the multiplicative Euler characteristic via a map from $K$-theory which makes sense without assuming Tate's semi-simplicity conjecture.  In particular, we can remove this hypothesis from a theorem of Milne proving a cohomological formula for zeta values attached to smooth proper $\fff_p$-schemes.  We also discuss extensions of these zeta value formulae to finite type $\fff_p$-schemes, and how recent progress in motivic homotopy theory allows us to prove some formulae without any assumptions on resolution of singularities or Tate's semi-simplicity conjecture.
\end{abstract}
{\small
	\setcounter{tocdepth}{1}
	\tableofcontents
	\vspace{3.0ex}
}
\section{Introduction}
\setcounter{section}{1}

\subsection{The Story so Far} Let $X$ be a scheme of finite type over a field $\fff_p$ of prime order $p$.  Recall that the Weil zeta function of $X$ is defined in \cite{weil1949numbers} as the unique meromorphic continuation to the complex plane of the product expansion $$\zeta(X,s)=\prod_{x\in X_{(0)}}\f{1}{1-|\kappa(x)|^{-s}},$$ where $X_{(0)}$ denotes the set of closed points of $X$, and $\kappa(x)$ is the residue field of a closed point, with $|\kappa(x)|$ denoting its order.  Following Grothendieck's solution to the Weil conjectures in \cite{grothendieck1964formule} (see also \cite{lubkin1968p} for the $p$-adic formulation), we learned that, at least for smooth proper $X$, the zeta function may be written as $\zeta(X,s)=Z(X,p^{-s})$, where \begin{equation}\label{Zetaboi}
Z(X,t)=\f{P_1(X,t)\ldots P_{2d-1}(X,t)}{P_0(X,t)\ldots P_{2d}(X,t)}, \qquad P_i(X,t)=\det(1-\varphi t|H^i(\overline{X},\bb{Q}_{\ell})).
\end{equation}  Here, $\varphi$ denotes the action of Frobenius on the $\ell$-adic cohomology of $X$, where we take $\ell$-adic \'{e}tale cohomology if $\ell\neq p$, and rational crystalline cohomology in the case $\ell=p$.

Following this, in \cite{Neukirch1978/79} Bayer and Neukirch were able to give, for smooth proper schemes $X/\fff_p$, formulae for $|\zeta(X,n)|_{\ell}$ with $\ell\neq p$ when $\zeta(X,s)$ does not have a pole at $n$.  As later noted by Schneider, their description follow essentially from linear algebra observations.  We will recall Schneider's argument in \textsection 3 as motivation for the arguments in the rest of the section.  This was extended by Schneider in \cite{Schneider1982} to describe, assuming a semi-simplicity conjecture due Tate, the value $$|\lim_{s\to n}(1-p^{n-s})^{\rho_n}\zeta(X,s)|_{\ell}$$ when $\zeta(X,s)$ has a pole at $n$ of order $\rho_n$.  Explicitly, the conjecture assumed is as follows, see for example the introduction to \cite{MilneValues} and also the conjecture $S(X,\ell)$ (really $S^r(X,\ell)$ for all $r$ in our formulation) from \cite{milne2007tate}.
\begin{conjecture}[Tate's semi-simplicity conjecture]\label{TateSemi}
For all primes $\ell\neq p$ (resp. $\ell=p$), all smooth proper schemes $X/\fff_p$, and all integers $n\in [0,\mo{dim}(X)]$, the Frobenius $\varphi$ acts semi-simply on the $\varphi=p^n$ eigenspace of $H^{2n}_{\acute{e}t}(\overline{X},\bb{Q}_{\ell})$ (resp. $H^{2n}_{\mo{crys}}(X,\bb{Q}_p)$ for $\ell=p$). 
\end{conjecture}  
We say that \textit{Tate's semi-simplicity conejecture holds for $X$} if the action of Frobenius on the $\ell$-adic cohomology of $X$ satisfies the conditions laid out in the conjecture above.  In the language of \cite{milne2007tate}, this is the assumption that $S^r(X,\ell)$ holds for all $r$.

In \cite{MilneValues}, Milne used p-adic Hodge theory, by studying the sheaves $W\Omega_{X,\mathrm{log}}^n[-n]$, to give a formula for computing the full value $|\zeta(X,n)|$ (or the respective leading coefficient), up to a sign, and the same semi-simplicity conjecture of Tate (to be completely precise, Milne only required the conjecture $S^n(X,\ell)$).
\begin{theorem}[\cite{MilneValues}, Theorem 0.1]\label{Milne1}
	Let $X/\fff_p$ be a smooth and proper scheme, and assume Tate's semi-simplicity conjecture holds for $X$.  Then for $n\in\zz$, we have:
	$$\zeta(X,s)\sim \pm \chi(X,\hat{\zz}(n),e)p^{\chi(X,\mc{O}_X,n)}(1-p^{n-s})^{-\rho_{n}}$$
	as $s\to n$.
\end{theorem}
\noindent Here, $H^i(X,\bb{Q}_{\ell})$ denotes the $\ell$-adic \'{e}tale cohomology of $\overline{X}$ for $p\neq \ell$ (treated as a $\mo{Gal}(\overline{\fff_p}/\fff_p)$-module), and the rational crystalline cohomology of $X$ when $p=\ell$.  We define $H^i_{\acute{e}t}(X,\hat{\zz}(n))=\prod_{\ell}H^i_{\acute{e}t}(X,\zz_{\ell}(n))$, where $H^i_{\acute{e}t}(X,\zz_{\ell}(n))$ denotes $\ell$-adic \'{e}tale cohomology of $X$ if $\ell\neq p$, and syntomic cohomology of $X$ if $\ell=p$.  The second term appearing in our formula, often called ``Milne's correcting factor,'' is defined as $p$ to the power of $\chi(X,\mc{O}_X,n)$, which itself is defined as
\begin{equation}\label{Milnecorrection}
\chi(X,\mc{O}_X,n):=\sum_{0\leq i\leq \dim(X), 0\leq j\leq n}(-1)^{i+j}(n-j)h^i(X,\Omega^j)
\end{equation} The first term, this $\chi(X,\widehat{Z}(n),e)$, which is well-defined whenever $n<0$ or $n>\mo{dim}(X)$, agrees with the multiplicative Euler characteristic of $R\Gamma(X,\widehat{\zz}(n))$ when this complex has finite cohomology groups.  When the complex doesn't have finite cohomology groups, under the assumption of semi-simplicity, the chain complex induced by the action of $e\in H^1(\fff_p,\widehat{\zz})$ on $R\Gamma(X,\widehat{\zz}(n))$, 
$$\ldots \to H^i(X,\widehat{\zz}(n))\xrightarrow{\cup e}H^{i+1}(X,\widehat{\zz}(n))\to\ldots$$ has finite cohomology groups, so the multiplicative Euler characteristic of this complex is well-defined, which Milne defines to be $\chi(X,\widehat{Z}(n),e)$ in \cite{MilneValues}.  We will recall this construction more formally in \textsection 2.1.

The formula from Theorem \ref{Milne1} was extended by Geisser in \cite{geisser2005arithmetic} to hold for all finite type schemes $X/\fff_p$ (of dimension $\leq d$), under the assumptions of a conjecture of Beilinson, Tate's semi-simplicity conjecture, and that there is a strong form of resolution of singularities in characteristic $p$ (up to dimension $d$),\footnote{Which we will recall in Definition \ref{GeisserRes}.} using an appropriate version of cohomology with compact support.  In the $p\neq \ell$ case, one could previously use compactly supported $\ell$-adic \'{e}tale cohomology to describe the $\ell$-adic absolute value of the special value at integer points, assuming Tate's semi-simplicity conjecture once again.  When one moves to the $\ell=p$ case, the theory of compactly supported p-adic cohomology theories is much less well-behaved, requiring the use of the so-called cdh topology which we motivate and recall in \textsection 4.2.  While we do not pursue studying the Weil-\'{e}tale topology version at present, one should be able to remove some assumptions from \cite{geisser2005arithmetic} by making use of the category we call $\cat{\zz}$.

\subsection{K-Theory to the Rescue} In order to get a formulation (and subsequent proof) of Theorem \ref{Milne1} which does not rely on Tate's semi-simplicity conjecture, we turn in \textsection 2.2 to defining, for any Dedekind ring $R$, a category $\cat{R}$ of arithmetic $C(S^1,R)$-modules.  We will only use this in the case $R=\zz_{\ell}$ here, but develop the general case as we believe it to be of independent interest, especially for Weil-\'{e}tale formulations of special value formulae and studying (possibly non-semi-simple) representation of the absolute Galois group $G_{\fff_p}$ on perfect $R$-modules for Dedekind rings $R$ (e.g. $R=\fff_q[T]$).  The key property of this category allowing us to define multiplicative Euler characteristics is summarized by the following theorem.
\begin{theorem}[\Cref{KZeroComp}]\label{KZeroCompI}
The 0th algebraic $K$-group of $\cat{R}$ splits as a direct sum of the 0th $K$-group of the subcategory of torsion modules and the free abelian subgroup generated by the unit.  That is, $$K_0(\cat{R})\simeq \zz\cdot [C(S^1,R)]\oplus K_0(\cat{R}_{\tors}).$$
\end{theorem}

Specializing to the case $R=\zz_{q}$ for a prime $q$ (either $\ell$ or $p$), we are able to define a version of this multiplicative Euler characteristic $\chi(-,e)$ (see also \Cref{IntegerEulerCharGen}).
\begin{definition}\label{EulerCharI}
We define the generalized multiplicative Euler characteristic $$\chi(-,e):K_0(\cat{\zz_{q}})\to \bb{Q}^{\times}$$ as the abelian group map determined by sending $[C(S^1,\zz_q)]\mapsto 1$ and factoring as the composite $K_0(\cat{\zz_q}_{\tors})\to K_0(\mo{Perf}(\zz_q)_{\tors})\xrightarrow{\chi(-)}\bb{Q}^{\times}$, with the first map the forgetful functor and the last map the classical multiplicative Euler characteristic.
\end{definition} 
Using this, we will prove, by way of a general recognition principal in Theorem \ref{DetectionTheorem}, that $\ell$-adic \'{e}tale cohomology and syntomic cohomology of smooth proper $\fff_{p}$-schemes canonically lift to $\cat{\zz_{q}}$ for $q$ either $\ell$ or $p$, and we are able to recover the previous versions of the multiplicative Euler characteristic $\chi(-,e)$ whenever they make sense.  In fact, we prove a bit more than this.

Let $\lge$ denote the dualizable objects in $\zz_{\ell}$-modules on the \'{e}tale site of $\fff_p$- that is, the bounded derived category of lisse $\zz_{\ell}$-sheaves on a point; and let $\fge=\fge(\fff_p)$ denote the category of perfect prismatic $F$-gauges on a point.  We prove the following.
\begin{theorem}[Theorems \ref{EtaleLifting} and \ref{FGaugeLift}, Proposition \ref{AgreementProp}]\label{LiftI}
There are canonical lifts of $$R\Gamma_{\acute{e}t}:\lge\to\mc{D}^{b}(\zz_{\ell}),\quad \text{ and } \quad R\Gamma_{\syn}:\fge\to \mc{D}^b(\zz_p)$$ to functors $$R\Gamma_{\acute{e}t}:\lge \to \cat{\zz_{\ell}} \quad \text{ and }\quad R\Gamma_{\syn}:\fge\to \cat{\zz_p}.$$
Moreover, given any complex $A\in \cat{\zz_{\ell}}$ such that Milne's multiplicative Euler characteristic $\chi^{\text{Milne}}(A,e)$ is well-defined, we have that $\chi^{\text{Milne}}(A,e)=\chi(A,e)$.
\end{theorem}
\begin{remark}\label{RemInfo}
These lifts carry slightly more information than the ordinary \'{e}tale cohomology would, but we argue that it is the right information to remember.  Indeed, in order to define Milne's version of the multiplicative Euler characteristic $\chi^{\text{Milne}}(A,e)$, you need to remember the action of cup product with $e$ on the cohomology of $A$.  For our theory, you roughly need to remember the action of $e$ itself on $A$ considered in the derived category, before passing to cohomology.
\end{remark}
\begin{remark}\label{RemMondal1}
A version of these multiplicative Euler characteristics $\chi(-,e)$ was defined independently by Mondal in \cite{mondal2025zetafunctionfgaugesspecial}, motivated by questions about prismatic F-gauges, via the so-called ``stable Bockstein characteristic.''  The stable Bockstein characteristic should provide an explicit model for computing our $\chi(-,e)$ in the case of Definition \ref{EulerCharI}, as long as one remembers the full $\zz$-action on a module rather than just the homotopy fixed points as a $C(S^1,\zz_{\ell})$-module (though in particular one should be able to take any lift of a $C(S^1,\zz_{\ell})$-module and use that to compute the stable Bockstein characteristic).  Mondal's results should apply more generally, with minor modifications, to provide a model for computing the version of the multiplicative Euler characteristic we define in Definition \ref{IntegerEulerCharGen}.
\end{remark}
\subsection{Special Value Formulae in the Smooth Proper Case}  With these results in hand, we proceed to show how one can use linear algebra arguments to determine the $q$-adic valuations of these special values of zeta functions, allowing us to provide a proof of Theorem \ref{Milne1} without requiring Tate's semi-simplicity conjecture.  To be precise, we prove the following.
\begin{theorem}\label{Main}
Let $X$ be a smooth proper $\fff_p$-scheme.  Then, for any integer $n$, we have that $$\zeta(X,s)\sim \pm \chi(X,\widehat{\zz}(n),e)p^{\chi(X,\mc{O}_X,n)}(1-p^{n-s})^{-\rho_n}$$ as $s\to n$.  Here $\chi(X,\widehat{\zz}(n),e)$ is defined to mean $\prod_{q\in\bb{P}}\chi(X,\zz_q(n),e)$ for primes $q$, with $\zz_{\ell}(n)$ denoting integral $\ell$-adic \'{e}tale cohomology, and $\zz_p(n)$ denoting syntomic cohomology.
\end{theorem}
\begin{proof}
This comes from combining Theorems \ref{ValuepMilne} and \ref{Valueell} in the case $n<0$ or $n>\mo{dim}(X)$, and Theorems \ref{MainTheoremEll} and \ref{MainTheoremValuep} in the case $0\leq n\leq \mo{dim}(X)$.
\end{proof}
\begin{remark}\label{RemMilne1}
We note that the above proof for $n<0$ or $n>\mo{dim}(X)$ is already due to Milne via Theorem \ref{Milne1}, and does not require our generalized version of the multiplicative Euler characteristic.  Also under the assumption of Tate's semi-simplicity conjecture (or even the weaker assumption $S^n(X,\ell)$), this simply recovers Milne's theorem \ref{Milne1}.
\end{remark}
\begin{remark}\label{RemMondal2}
The results above in the $p$-completed case, specifically Theorem \ref{MainTheoremValuep}, have appeared previously in work of Mondal \cite{mondal2025zetafunctionfgaugesspecial} under a different guise, using a stable Bockstein characteristic defined in loccit as opposed to our generalization $\chi(-,\zz_p(n),e)$.  Although Mondal does not pursue it in \cite{mondal2025zetafunctionfgaugesspecial}, one should also be able to recover a variant of Theorem \ref{MainTheoremEll} using his stable Bockstein characteristic in place of our generalized multiplicative Euler characteristic. 
\end{remark}
\begin{remark}
The proof of Theorem \ref{Main} applies not only for zeta functions attached to smooth proper $\fff_p$-schemes, but for those attached to some suitably nice category of ``motives'' on $\mo{Spec}(\fff_p)$, at least defined such that $M_{\ell}\in \lge$ for primes $\ell\neq p$, and $M_p\in \fge$ in some natural way.  By virtue of having a generalized Euler characteristic defined for arithmetic $C(S^1,R)$-modules over any Dedekind ring $R$, this should allow one to even perform similar techniques integrally, and perhaps even for generalizations to ``motives with coefficients in nice enough rings.''
\end{remark}
\subsection{Special Value Formulae for Finite Type $\fff_p$-Schemes}  In \textsection 4, we move away from the smooth proper case, and focus on proving special value formulae for arbitrary finite type $\fff_p$-schemes.  When $\ell\neq p$, the analogue of Theorem \ref{MainTheoremEll} works on the nose, using compactly supported \'{e}tale cohomology, and we can prove, without assuming Tate's semi-simplicity conjecture,
\begin{theorem}[Theorem \ref{GeneralizationFTell}]\label{FTellI}
Let $X$ be any finite type $\fff_p$-scheme.  Then, $R\Gamma_{\acute{e}t}(X,\zz_{\ell})_{c}$ has a canonical lift to $\cat{\zz_{\ell}}$, and we can compute $$|C(X,n)|_{\ell}^{-1}=\chi(R\Gamma(X,\zz_{\ell}(n))_c,e),$$
where $C(X,n)$ is the rational number such that $\zeta(X,s)\sim (1-p^{n-s})^{-\rho_n}C(X,n)$ as $s\to n$, for $\rho_n$ the order of the pole of $\zeta(X,s)$ at $s=n$.
\end{theorem}
\begin{remark}\label{RemTate2}
In the presence of Tate's semi-simplicity conjecture, this was studied (in a slightly different guise) by Geisser in \cite{geisser2005arithmetic}, and under this assumption, our result reduces to the ``$\ell$-adic completion of his result.''  Geisser worked with the Weil-\'{e}tale topology, studying all primes at once, whereas we do it one prime at a time.  The correct multiplicative Euler characteristic to use to recover Geisser's results without Tate's semi-simplicity conjecture would be the integral variant defined in Definition \ref{IntegerEulerCharGen}.
\end{remark}
When we move to the case $\ell=p$, we do not get nice compactly supported syntomic/Hodge cohomology groups, and it's not even immediately obvious how to define such compactly supported cohomology groups in the first place.  We will use this fact to motivate the \textit{cdh topology}, defined originally by Suslin and Voevodsky in \cite{suslin2000relative}, and it provides a way to discuss ``compactly supported cohomology.''  Under an assumption of strong resolution of singularities (in dimension $\leq d$), one can prove special value formulae for general finite type $\fff_p$-schemes (of dimension $\leq d$) using cdh descent.  The best we can do absent resolution of singularities, using only cdh sheafifications, is summarized by the following theorem.

\begin{theorem}\label{Main2}
	Let $\mc{G}\subseteq \mo{Sch}^{prop}_{\fff_p}$ be the smallest subcategory of proper $\fff_p$-schemes which contains all smooth proper $\fff_p$-schemes $X$ for which the map induced by sheafification
	\begin{equation}
	 \mo{R}\Gamma(X,W\Omega/\mc{N}^{\geq k}W\Omega)\to L_{cdh}\mo{R}\Gamma(X,W\Omega/\mc{N}^{\geq k}W\Omega),
	\end{equation}
	is an equivalence, and which has the property that given an abstract blowup square 
	\begin{center}
		\begin{tikzcd}
			Z^{\prime}\dar\rar & X^{\prime}\dar\\
			Z\rar& X
		\end{tikzcd}
	\end{center}
	with any three of the displayed schemes in $\mc{G}$, then the fourth is as well.
	
	Then, for any finite type $\fff_p$-scheme $U$ having a compactification $U\to X$, with closed complement $Z$, such that $X,Z\in \mc{G}$, $\chi((L_{cdh}\zz_p(n))_c(U),e)$ is well-defined, and we have 
	\begin{equation}
		|C(U,n)|_p^{-1}=\chi((L_{cdh}\zz_p(n))_c(U),e)\chi((L_{cdh}W\Omega_{-}/\mc{N}^{\geq n}W\Omega_{-})_c(U)).
	\end{equation}
\end{theorem}
\begin{proof}
This follows from Theorem \ref{MainThFT}, using work of Elmanto-Morrow \cite[Corollary~6.5]{elmanto2023motivic}, proving that the natural map $$R\Gamma_{\syn}(X,\zz_p(n))\to L_{cdh}R\Gamma_{\syn}(X,\zz_p(n))$$ is an equivalence for every smooth proper $\fff_p$-scheme $X$.
\end{proof}
\begin{remark}\label{RemGeisser1}
Under the assumption of a strong form of resolution of singularities in characteristic $p$, and Tate's semi-simplicity conjecture, the class $\mc{G}$ contains all proper $\fff_p$-schemes, and the result above effectively reduces to a result of Geisser in \cite{geisser2005arithmetic}.  In \textsection 4.2, we will indicate in more detail how our results relate to those of Geisser (see Corollary \ref{GeisserTheoremCor}).
\end{remark}
The strongest possible result we could prove absent full resolution of singularities would be if we knew that the maps $$\mo{R}\Gamma(X,W\Omega/\mc{N}^{\geq k}W\Omega)\to L_{cdh}\mo{R}\Gamma(X,W\Omega/\mc{N}^{\geq k}W\Omega)$$ were equivalences for every smooth proper $\fff_p$-scheme $X$ and every $k$.  This would follow for instance (and is in fact equivalent to) an analogous claim for Hodge cohomology- that $$\mo{R}\Gamma(X,\Omega^j)\to L_{cdh}\mo{R}\Gamma(X,\Omega^j)$$ is an equivalence for every smooth proper $\fff_p$-scheme $X$.  Similar results before taking derived functors have been proven by Huber-Kelly \cite{Huber_Klawitter_2018} (see also Ertl-Miller \cite{ERTL20195285}).

Nevertheless, we can still improve upon Theorem \ref{Main2} absent resolution of singularities.  Namely, recent work of Annala-Pstragowski \cite{annala2025noteweightfiltrationscharacteristic} and Annala-Hoyois-Iwasa \cite{annala2024atiyahdualitymotivicspectra} constructs a motivic spectrum representing $R\Gamma(-,\Omega^j)$, whose $\bb{A}^1$-colocalization, which we denote by the dagger notation $R\Gamma(-,(\Omega^j)^{\dagger})$, agrees with the original spectrum on smooth \textit{projective} $\fff_p$-schemes.  Since $\bb{A}^1$-invariant motivic spectra determine cdh-local cohomology theories by \cite{Khan_2024} and \cite{Cisinski_2013}, we can prove the following theorem without any assumption on resolution of singularities or Tate's semi-simplicity conjecture.
\begin{theorem}[Theorem \ref{Strongestuncondtheorem}]\label{ToniExtension}
Let $\mc{G}^{\prime}\subseteq \mo{Sch}^{\mathrm{prop}}_{\fff_p}$ be the smallest full subcategory of proper $\fff_p$-schemes which contains all smooth proper $\fff_p$-schemes $X$ such that $$R\Gamma(X,(\Omega^j)^{\dagger}) \to R\Gamma(X,\Omega^j)$$ is an equivalence for all $j$, and which is closed under the three-out-of-four for abstract blowup squares (\ref{absblowup}).  Then for any finite type $\fff_p$-scheme $U$ which has a compactification $U\to X$ with closed complement $Z$ such that $X,Z\in\mc{G}^{\prime}$, we have 
$$
|C(U,n)|_p^{-1}=\chi((L_{cdh}\zz_p(n))_c(U),e)\cdot p^{\chi(U,\mc{O}_U,n)_c^{\dagger}},
$$
where $$\chi(U,\mc{O}_U,n)_c^{\dagger}:=\sum_{j=0}^{n}(-1)^{j}(n-j)\chi(R\Gamma(U,(\Omega^j)^{\dagger})_c).$$
In particular, if we can find a smooth projective compactification of $U$ with smooth projective complement, this formula applies.
\end{theorem}

\subsection{Conventions and Notation}  We freely use the theory of (stable) $\infty$-categories in the sense of Lurie \cite{HA} and \cite{lurie2008highertopostheory}.  In particular, all derived categories will be assumed to be derived $\infty$-categories, so that we can take co/limits, and module categories over $\bb{E}_{\infty}$-rings will be taken to mean module spectra over said rings.  When discussing certain sheaves on schemes $X$ such as the K\"{a}hler differentials $\Omega^j_{X/\fff_p}$, we will often abusively identify $\Omega^j_{X/\fff_p}$ with $R\Gamma(X,\Omega^j_{X/\fff_p})\in\mc{D}^b(\fff_p)$, since we will only ever really deal with cohomology in this paper.

\subsection{Acknowledgments} I would like to thank Baptiste Morin and Don Blasius for encouraging me to type up the original version of this note, and for helpful conversations.  I would also like to thank Brian Shin, Matthew Morrow, and Jas Singh for helpful conversations related to this work, and Elden Elmanto for comments on the previous versions.  I would like to thank Toni Annala for telling me about the $\bb{A}^1$-colocalization procedure which lead to the extension of Theorem \ref{MainThFT} to Theorem \ref{Strongestuncondtheorem}.  Finally, I would like to thank Adrien Morin, Baptiste Morin, and Oakley Edens for helpful comments on the current version of the paper, and for discussions which were of great help with writing up the exposition for this work.

I would like to especially thank Shubhodip Mondal for pointing out a glaring error in a previous version's definition of Arithmetic $C(S^1,R)$-modules, for which the analogue of \Cref{DetectionTheorem} was false, which this version quickly corrects.  In particular, this also modifies the result \Cref{KZeroComp}, although the main applications remain unaffected.

\newpage
\section{Multiplicative Euler Characteristics}
\indent\indent In this section, we recall some standard homological algebra arguments used throughout this paper, simultaneously introducing some notational conventions we choose to take.  A major player in this paper is the multiplicative Euler characteristic for certain complexes. 
\subsection{Classical Multiplicative Euler Characteristics}  We begin by recalling the classical story of multiplicative Euler characteristics.  We discuss the story only for complexes defined over the integers $\zz$, but with the obvious modifications one can make sense of this with the $p$-adic integers $\zz_p$ replacing $\zz$ everywhere, which will be what we primarily use in this paper.
\begin{definition}\label{ClassicDef1}
Let $A\in\mc{D}^b(\zz)$ be a complex the bounded derived category of $\zz$, such that the cohomology group $H^i(A)$ is a finite abelian group for all $i\in\zz$.  The \textit{multiplicative Euler characteristic} $\chi(A)$ is defined by $$\chi(A)=\prod_{i\in\zz}|H^i(A)|^{(-1)^i},$$ where $|H^i(A)|$ denotes the order of the (finite) group $H^i(A)$.
\end{definition}
Since we assumed $A$ was in the bounded derived category of $\zz$, $H^i(A)=0$ for $|i|\gg 0$, so that this multiplicative Euler characteristic is well-defined.  We record an immediate consequences of this definition,
\begin{lemma}\label{TrivialClaim1}
Let $A\to B \to C$ be a fiber sequence in $\mc{D}^b(\zz)$ with $A,B,C$ all having finite cohomology groups, such that their multiplicative Euler characteristics are defined.  Then we have that $\chi(B)=\chi(A)\chi(C)$.
\end{lemma}
\begin{proof}
This follows from the long exact cohomology sequence and counting cardinalities of images/kernels.
\end{proof}
One fact we will have to make use of in order to relate our correction factor to Milne's is how the multiplicative Euler characteristic interacts with finite filtrations.  Recall the following definition from \cite[\textsection Appendix B]{burklund2022galois} and \cite[\textsection 1.2.2]{HA}, to which we refer for further details.
\begin{definition}\label{Filteredbois}
A \textit{filtered object} in $\mc{D}(\zz)$ (resp. a general stable $\infty$-category $\mc{C}$) is a functor $$F^*:\zz^{\leq}\to \mc{D}(\zz)$$ (resp. $F^*:\zz^{\leq}\to\mc{C}$) from the integers considered as a poset to the derived $\infty$-category $\mc{D}(\zz)$ (resp. $\mc{C}$).  We say that $F^*$ is a \textit{finite filtration} on a complex $A\in \mc{D}(\zz)$ if $F^n=0$ for $n\ll 0$ sufficiently small, $F^n\simeq A$ for $n\gg0$, and the transition maps $F^n\to F^{n+1}$ are equivalences for all $n \gg 0$.
\end{definition}
\begin{definition}\label{taumoment}
Given a filtered object $F^*:\zz^{\leq}\to\mc{D}(\zz)$, the \textit{$n$th associated graded} piece of $F^*$, denoted $\mo{gr}^n(F)$, is given by the cofiber $\mo{cofib}(F^{n-1}\to F^n)$.
\end{definition}
Now we come to the main proposition that will be used to compare our correction factor to Milne's.
\begin{proposition}\label{trivialclaim2}
Suppose $A\in \mc{D}^b(\zz)$ is equipped with a finite filtration $F^*:\zz^{\leq }\to \mc{D}^b(\zz)$ such that for all $n\in\zz$, $F^n\in\mc{D}^b(\zz)$ has finite cohomology groups, so that multiplicative Euler characteristics can be defined for each of these complexes.  Then the same is true of $\mo{gr}^n(F)$, and $$\chi(A)=\prod_{n\in\zz}\chi(\mo{gr}^n(F)).$$
\end{proposition}
\begin{proof}
By our assumption that $F$ is a finite filtration, $\mo{gr}^n(F)$ is zero for almost all $n$, so that the product on the right is well-defined.  The result now follows from Lemma \ref{TrivialClaim1} and induction, using the fiber sequence $F^{n-1}\to F^n\to \mo{gr}^n(F)$ to write $$\chi(F^n)=\chi(F^{n-1})\chi(\mo{gr}^n(F)),$$ where we use that $F^n=0$ for $n\ll0$ to begin the induction.
\end{proof}

Now, there are two major pieces of notation from the introduction that we have yet to explain.  If we have a presheaf $\mc{F}:\mo{Sch}^{ft,op}_{\fff_p}\to\mc{D}^b(\zz)$ which on a particular scheme $X$ is such that $\mc{F}(X)$ has finite cohomology groups, we will sometimes write $\chi(X,\mc{F})$ to mean $\chi(\mc{F}(X))$.  This leaves only the notation $\chi(X,\zz_p(n),e)=\chi(\zz_p(n)(X),e)$ from the statement of Theorem \ref{Milne1}.  We define this now, as it will be important for describing special values at poles of $\zeta(X,s)$.
\begin{definition}\label{classicaladhocdef}
Suppose that $A\in\mc{D}^b(\zz)$ is an arbitrary complex, equipped with a map $e:A\to \Sigma^{-1} A$ with $e^2\simeq 0$.  Assume further that the cohomology groups of the complex $(H^*(A),e)$ induced from $e$: $$\ldots\xrightarrow{e} H^i(A)\xrightarrow{e} H^{i+1}(A)\xrightarrow{e} \ldots$$ are all finite.  Then we define the \textit{multiplicative Euler characteristic of the pair $(A,e)$} to be the multiplicative Euler characteristic of the complex $(H^*(A),e)$ built from the action of $e$, that is $$\chi(A,e):=\prod_{j\in\zz}|H^j((H^*(A),e))|^{(-1)^j}.$$
\end{definition}
We make the following observation.
\begin{lemma}\label{classicalagreementlemma}
If $A\in\mc{D}^b(\zz)$ has finite cohomology groups, then for any map $e:A\to\Sigma^{-1}A$ with $e^2\simeq 0$, $\chi((A,e))=\chi(A)$.
\end{lemma}
\begin{proof}
The $i$th cohomology group of the complex $(H^*(A),e)$ is given by $$\mo{ker}(e:H^i(A)\to H^{i+1}(A))/\mo{im}(e:H^{i-1}(A)\to H^i(A)).$$  Since these are all finite groups, we have that $$|H^j((H^*(A),e))|=|\mo{ker}(e:H^i(A)\to H^{i+1}(A))|\cdot |\mo{im}(e:H^{i-1}(A)\to H^i(A))|^{-1}.$$  Furthermore, again since $H^i(A)$ is a finite group, we have that $$|H^i(A)|=|\mo{ker}(e:H^i(A)\to H^{i+1}(A))|\cdot |\mo{im}(e:H^{i}(A)\to H^{i+1}(A))|.$$  Putting these together, we find that 
\begin{align*}
	\chi(A)&=\prod_{i\in\zz}|H^i(A)|^{(-1)^{i}}\\
	&=\prod_{i\in\zz}|\mo{ker}(e\colon H^i(A)\to H^{i+1}(A))|^{(-1)^{i}}\cdot |\mo{im}(e\colon H^{i}(A)\to H^{i+1}(A))|^{(-1)^{i}}\\
	&=\prod_{i\in\zz}|\mo{ker}(e\colon H^i(A)\to H^{i+1}(A))|^{(-1)^{i}}\cdot |\mo{im}(e\colon H^{i-1}(A)\to H^{i}(A))|^{(-1)^{i-1}}\\
	&=\prod_{i\in\zz}\left(|\mo{ker}(e\colon H^i(A)\to H^{i+1}(A))|\cdot |\mo{im}(e\colon H^{i-1}(A)\to H^{i}(A))|^{-1}\right)^{(-1)^{i}}\\
	&=\prod_{i\in\zz}|H^i((H^*(A),e))|^{(-1)^i}=\chi(A,e),
\end{align*}
as claimed.
\end{proof}
\subsection{A K-theoretic Variant}  Our goal now is to redefine a version of $\chi(A,e)$ that will make sense for certain complexes $A\in\mc{D}^b(\zz)$ equipped with a square-zero map $e:A\to\Sigma A$, even when the definition 2.6 does not necessarily make sense.  For the complexes we consider in the sequel, the map $e:A\to \Sigma A$ arises from a module structure for $A$ over the $\bb{E}_{\infty}-\zz-$algebra $\zz^{h\zz}=C(S^1,\zz)$,\footnote{Or rather, in our cases, $C(S^1,\zz_p)$, but again we proceed with the integral construction noting that it works exactly the same in the $p$-complete world.} which has $H^i(C(S^1,\zz))=\zz$ if $i=0,1$ and $0$ otherwise.  What we would like to say is that $\chi(-,e)$ is additive along fiber sequences whenever it is defined, and it satisfies the 2-out-of-3 property for being defined on $C(S^1,\zz)$-modules.  Unfortunately, this claim is simply false, as we have the cofiber sequence $$\Sigma^{-1}C(S^1,\zz)\to  C(S^1,\zz)\to N_2,$$ using notation from \cite[\textsection 3.3]{mathew2016residuefieldsclassrational}, where $H^i(N_2)=\zz$ if $i=0,1$, and is zero otherwise, but where naturality forces $e$ to act as $0$ on $H^i(N_2)$.  Explicitly, $N_2$ arises as the homotopy fixed points $V^{h\zz}$, where $V\simeq \zz\oplus \zz$, and $\zz$ acts on $V$ through the matrix $\begin{bmatrix}
	1 & 1 \\ 0 & 1
\end{bmatrix}$, that is to say, $N_2$ witnesses a failure of semi-simplicity for the $\zz$ action on $V$.

Thus, we are forced to consider an alternative.  To motivate our following definition, we reinterpret the multiplicative Euler characteristic from definition \ref{ClassicDef1}.  Recall that a bounded complex $A\in \mc{D}^b(\zz)$ has finite cohomology groups if and only if $A$ lives in the full subcategory $\mo{Perf}(\zz)_{\mathrm{tors}}$ of $\mo{Perf}(\zz)$, which arises as the kernel of the localization functor $\mo{Perf}(\zz)\to\mo{Perf}(\bb{Q})$.  We note the following.
\begin{lemma}\label{Kzeroequivalence1}
There is an equivalence on $K$-theory groups $K_0(\mo{Perf}(\zz)_{\mathrm{tors}})\simeq \bigoplus_{p\in\bb{P}}\zz$, where $\bb{P}$ denotes the set of primes, with $\fff_p$ generating the $p$-indexed copy of $\zz$.  Moreover, there is a unique group homomorphism $\chi^K:K_0(\mo{Perf}(\zz)_{\mathrm{tors}})\to \bb{Q}^{\times}$ such that given a complex $A\in \mo{Perf}(\zz)_{\mathrm{tors}}$, representing the class $[A]$ in $K_0$, $\chi(A)=\chi^K([A])$.
\end{lemma}
\begin{proof}
The first claim is standard.  To define $\chi^K$, it suffices to define it on free generators of $K_0$, which we do by stipulating that $\chi^K(\fff_p)=p$.  By using the structure theorem for finite abelian groups, the fiber sequences $\zz/p^{n-1}\zz\to \zz/p^n\zz\to\zz/p\zz$, and that any complex of $\zz$-modules is quasi-isomorphic to its cohomology, we conclude that for $A\in \mo{Perf}(\zz)_{\mathrm{tors}}$, $\chi(A)=\chi^K([A])$.
\end{proof}
This lemma allows us to justify conflating the maps $\chi$ and $\chi^K$.  

When trying to adapt this to multiplicative Euler characteristics attached to complexes with a map, we almost want to replicate the above lemma with $\mo{Perf}(C(S^1,\zz))$ replacing $\mo{Perf}(\zz)_{\mathrm{tors}}$.  This however, does not quite work, since once can use work of Burklund-Levy \cite{burklund2021ktheoryregularcoconnectiverings} to show that $K_0(\mo{Perf}(C(S^1,\zz)))=K_0(\zz)=\zz$, generated by the unit.  Luckily, this category is also not the correct one to look at, since, if we take the $\ell$-completed analogue, $R\Gamma_{\acute{e}t}(X,\zz_{\ell}(n))$ is very rarely perfect over $C(S^1,\zz_{\ell})$ even for $X$ smooth proper.  In fact, going so far as to even take $X=\mo{Spec}(\fff_p)$ to be a point, this \'{e}tale cohomology will fail to be perfect over $C(S^1,\zz_{\ell})$ for general $n$.  Thus, we enlarge the category $\mo{Perf}(C(S^1,\zz))$ slightly in order to define our generalized Euler characteristics.

Since the following constructions hold in suitable generality, we will state them for an arbitrary Dedekind domain $R$ in place of $\zz$ (or $\zz_{\ell}$), specializing to these cases later when defining our generalization of the multiplicative Euler Characteristic.
\begin{definition}\label{arithmeticmods}
Let $R$ be a Dedekind domain.  The category of \textit{arithmetic $C(S^1,R)$-modules} is defined as the smallest full thick stable subcategory $\cat{R}\subseteq \mo{Mod}(C(S^1,R))$ of $C(S^1,R)$-modules which contains both $\mo{Perf}(C(S^1,R))$ and all modules $M$ over $C(S^1,R)$ whose underlying $R$-module is both perfect and torsion.
\end{definition}
Let's define some explicit perfect modules over $C(S^1,R)$, following the discussion in \cite[Construction 3.19]{mathew2016residuefieldsclassrational}, from which we borrow the notation.  The following modules are meant to capture the information of ``Jordan blocks.''
\begin{definition}[\cite{mathew2016residuefieldsclassrational} Construction 3.19]\label{AkhilDef}
For $m\geq 1$, let $N_m$ denote the $C(S^1,R)$-module obtained by taking the homotopy fixed points $(V_m)^{h\zz}$ for $V_m=R[x]/x^m$ in homological degree 0, where a generator $\sigma\in\zz$ acts through the automorphism $1+x$ of $V_m$.
\end{definition}
We now prove these modules actually are in fact perfect, and sit in rather convenient fiber sequences.
\begin{lemma}\label{LemAkhilDef}
The module $N_m$ has $$\pi_i(N_m)=\begin{cases}
	R & \qquad \text{ if }i=0,-1\\
	0 & \qquad \text{ otherwise.}
\end{cases}$$  Moreover, taking the map $\Sigma^{-1}C(S^1,R)\to N_m$ picking out a generator of $\pi_{-1}(N_m)$, we have a cofiber sequence $$\Sigma^{-1}C(S^1,R)\to N_m\to N_{m+1}.$$ In particular, since $N_1\simeq C(S^1,R)$, the modules $N_m$ are perfect for all $m\geq 1$.
\end{lemma}
\begin{proof}
It suffices to show the existence of the claimed cofiber sequences.  Note that, for $m\geq 1$, there is a $\zz$-equivariant fiber sequence $$V_m\xrightarrow{x} V_{m+1}\to V_1,$$ where $V_1^{h\zz}\simeq C(S^1,R)$.  Since taking homotopy fixed points is exact, this leads to a fiber sequence $$N_m\to N_{m+1}\to C(S^1,R),$$ and the long exact sequence in homotopy groups induced from this shows that the map $\Sigma^{-1}C(S^1,R)\to N_m$ induced from this fiber sequence picks out a generator of $\pi_{-1}(N_m)$, as claimed.
\end{proof}
\begin{corollary}\label{NiceCorCircle}
The direct limit of the system defined by the maps in the previous lemma is $R$, considered with the trivial action, that is $$\varinjlim(\ldots\to N_m\to N_{m+1}\to\ldots)\simeq R.$$
\end{corollary}
\begin{proof}
This can be seen by using commutation of homotopy groups with filtered colimits, but we provide another argument as well.  The functor $(-)^{h\zz}$ is a finite limit, and in particular commutes with colimits since we work in a stable $\infty$-category.  Thus, $$\varinjlim N_m\simeq (\varinjlim V_m)^{h\zz}\simeq (R[x^{\pm 1}]/xR[x])^{h\zz}\simeq R.$$
\end{proof}
\begin{remark}
While the above can be used to provide a unique $C(S^1,R)$-module structure on the $R$-module $R$, this does \textit{not} imply the existence of unique $C(S^1,R)$-module structures on modules such as $R/a$, as was claimed in a previous version of this paper (which was kindly pointed out be false to us by Mondal).  The reason for this is that $\mo{End}_{C(S^1,R)}(R)\simeq R[[x]]$, and the forgetful functor sends $x\mapsto 0$, so one can take quotient e.g. the cofiber by $a-x\colon R\to R$ which has underlying $R$-module $R/a$ but is \textit{not} $R/a$ with the trivial action as a $C(S^1,R)$-module.
\end{remark}
Since $\cat{R}$ is defined a little abstractly, we record the following key criterion for checking that a module $A\in\mo{Mod}(C(S^1,R))$ lives in $\cat{R}$.
\begin{theorem}\label{DetectionTheorem}
Let $A\in \mo{Mod}(C(S^1,R))$ be a module over the $\bb{E}_{\infty}$-ring $C(S^1,R)$.  Then we have that $A\in \cat{R}$ if and only if $\pi_*(A)$ is a finitely generated (graded) $R$-module and $A\otimes_{R}F$ is perfect over $C(S^1,F)$, where $F$ is the fraction field of $R$.
\end{theorem}
\begin{proof}
Since localizing from $R$ to $F$ commutes with the forgetful functor from $C(S^1,R)$-modules back to $R$-modules, the rationalization of any torsion arithmetic module is zero.  In particular, for any $M\in\cat{R}$, then $M\otimes_{R}F$ is in the stable subcategory generated by perfect $C(S^1,F)$-modules and $0$, which is to say that it is perfect.  All of the generators of $\cat{R}$ have perfect underlying $R$-module, so the same is true of any module in $\cat{R}$, telling us that the cohomology groups of any such $M$ are finitely generated, with only finitely many nonzero, giving the first description.

Conversely, suppose we are given some $A\in\mo{Mod}(C(S^1,R))$ with finitely generated cohomology groups such that $A\otimes_{R} F$ is perfect.  Choose some perfect $C(S^1,R)$-module $P$ with $P\otimes_{R}F\simeq A\otimes_{R}F$, which we can construct explicitly using e.g. the classification of perfect $C(S^1,F)$-modules from \cite{mathew2016residuefieldsclassrational}.  Fixing an equivalence $$N\otimes_{R}F\simto A\otimes_{R}F,$$ and composing with $N\to N\otimes_{R}F$, we obtain a map $$N\to A\otimes_{R}F$$ which becomes an equivalence after rationalization.  Since $$A\otimes_{R}F\simeq \colim_{f_1,\ldots,f_n\in R\bs \{0\}}A\xrightarrow{f_1\ldots f_n}A$$ is a filtered colimit along multiplication maps witnessing this localization, and $N$ is a compact object in $\mo{Mod}(C(S^1,R))$, we find that $N\to A\otimes_{R}F$ lifts to a map $N\to A$ which becomes an equivalence after rationalization.  In particular, we can form $$B:=\mo{cofib}(N\to A),$$ which is the cofiber of a map from a perfect module to $A$, so in particular is in $\cat{R}$ if and only if $A$ is.

Since $B$ rationalizes to zero, and has homotopy groups finitely generated over $R$, we find that $\pi_*(B)$ is a finitely generated torsion $R$-module.  Since $R$ is a Dedekind ring, it has global dimension 1, and in particular any $R$-module with finitely generated cohomology groups (which are zero in all but finitely many degrees) is perfect, giving that $B$ has underlying $R$-module which is both torsion and perfect, and therefore lies in $\cat{R}$.

\end{proof}

Roughly speaking, our goal with this definition is to include both complexes where ``some variant of $e$'' can be used to remove the torsion-free part (in particular we do not want to include the $C(S^1,\zz)$-module $\zz$), but also allow for torsion modules which we can measure the size of.

Now, we prove the main theorem which will allow us to define ``$\chi(A,e)$'' for arbitrary $A\in \cat{\zz}$, and similarly for $q$-completed variants.
\begin{theorem}\label{KZeroComp}
We have a canonical splitting of Grothendieck groups $$K_0(\cat{R})\simeq K_0(\cat{R}_{\tors})\oplus \zz\cdot[\bb{1}].$$
\end{theorem}
\begin{proof}
Note that the localization sequence $$\cat{R}_{\mathrm{tors}}\to \cat{R}\to \mo{Perf}(C(S^1,F))$$ induces an exact sequence $$K_1(C(S^1,F))\to K_0(\cat{R}_{\mathrm{tors}})\to K_0(\cat{R})\to K_0(C(S^1,F)).$$  By \cite[Theorem 1.1]{burklund2021ktheoryregularcoconnectiverings}, $K_i(C(S^1,F))\simeq K_i(F)$ for $i\geq 0$, which follows from the facts that $C(S^1,F)$ is a coconnective $\bb{E}_1$-ring (which is even $\bb{E}_{\infty}$), $\pi_0(C(S^1,F))=F$ is a field (and is in particular regular coherent), and $\pi_{-1}(C(S^1,F))$ is of course flat over this field.  Since $K_0(C(S^1,F))\simeq \bb{Z}$, generated by the class of the unit $[\bb{1}]$, the map sending $[\bb{1}]\mapsto[\bb{1}]$ splits the surjection $K_0(\cat{R})\to K_0(C(S^1,F))$.  To prove our desired claim, it suffices to identify the kernel of this surjection.

Towards this end, we examine the following commutative diagram of stable $\infty$-categories
\begin{center}
	\begin{tikzcd}
		\mo{Perf}(C(S^1,R))_{\mathrm{tors}}\rar\dar & \mo{Perf}(C(S^1,R))\rar\dar & \mo{Perf}(C(S^1,F))\ar[d,equal]\\
		\cat{R}_{\mathrm{tors}}\rar & \cat{R}\rar & \mo{Perf}(C(S^1,F)).
	\end{tikzcd}
\end{center}
This induces the following commutative diagram of $K$-groups, with exact rows
\begin{center}
	\begin{tikzcd}
		K_1(C(S^1,F))\ar[d,equal]\rar &K_0(\mo{Perf}(C(S^1,R))_{\mathrm{tors}})\rar\dar & K_0(C(S^1,R))\rar\dar & K_0(C(S^1,F))\ar[d,equal]\\
		K_1(C(S^1,F))\rar & K_0(\cat{R}_{\mathrm{tors}})\rar & K_0(\cat{R})\rar & K_0(C(S^1,F)).
	\end{tikzcd}
\end{center}
The claim we are after reduces now to showing that the map $K_1(C(S^1,F))\to K_0(\cat{R}_{\tors})$ is the zero map.  By commutativity of the above diagram, it suffices to show that the map $K_0(\mo{Perf}(C(S^1,R))_{\tors})\to K_0(\cat{R}_{\tors})$ is zero.

The category $\mo{Perf}(C(S^1,R))_{\mathrm{tors}}$ has orthogonal generators $C(S^1,R/\mf{p})$ for height 1 primes $\mf{p}$ of $R$, with $$\pi_i(\mo{End}_{\mo{Perf}(C(S^1,R))_{\mathrm{tors}}}(C(S^1,R/\mf{p})))=\begin{cases}
	R/\mf{p} & \qquad \text{ if }i=0,-2\\
	R/\mf{p}\oplus R/\mf{p} & \qquad \text{ if } i=-1\\
	0 & \qquad \text{ otherwise.}
\end{cases}$$  This is a coconnective $\bb{E}_1$-ring, with $\pi_0$ given by a field, and $\pi_{-i}$ necessarily flat over it, so \cite[Theorem 1.1]{burklund2021ktheoryregularcoconnectiverings} applies.  Using a Schwede-Shipley argument, we find that $$K_0(\mo{Perf}(C(S^1,R))_{\mathrm{tors}})\simeq \bigoplus_{\mf{p}\in\mo{Spec}(R)_{(0)}} \zz,$$ with the $\mf{p}$-indexed copy of $\zz$ generated by $\left(R/\mf{p}\right)^{h\zz}=C(S^1,R/\mf{p})$.

We have an explicit description of the map $K_0(\mo{Perf}(C(S^1,R))_{\mathrm{tors}})\to K_0(\cat{R}_{\tors})$, as it takes the free generator $\left(R/\mf{p}\right)^{h\zz}$ to its class $[\left(R/\mf{p}\right)^{h\zz}]\in K_0(\cat{R}_{\mathrm{tors}})$.  To compute what this class is, we simply use the fiber sequence $$\Sigma^{-1}R/\mf{p}\to \left(R/\mf{p}\right)^{h\zz}\to R/\mf{p},$$ which tells us that in $K_0$, one has $$[\left(R/\mf{p}\right)^{h\zz}]=[R/\mf{p}]+[\Sigma^{-1}R/\mf{p}]=[R/\mf{p}]-[R/\mf{p}]=0,$$ as desired.

\end{proof}

Specializing now to the case $R=\zz$, this leads us to define
\begin{definition}\label{IntegerEulerCharGen}
Let $\chi(-,e):K_0(\cat{\zz})\to\bb{Q}^{\times}$ be the group homomorphism which factors as the composite $$K_0(\cat{\zz})\to K_0(\cat{\zz}_{\tors})\to K_0(\mo{Perf}(\zz)_{\tors})\xrightarrow{\chi}\bb{Q}^{\times},$$
where the first map is the projection induced from \Cref{KZeroComp} sending $[\bb{1}]$ to zero, the second map is induced by the forgetful functor, and the last map is the usual multiplicative Euler characteristic from \Cref{ClassicDef1}.  For any complex $A\in\cat{\zz}$, we often write $\chi(A,e)$ for $\chi([A],e)$.
\end{definition}
Finally, we show that given $A\in\cat{\zz}$ such that the induced action of $e\in \pi_{-1}(C(S^1,\zz))$ on $A$ satisfies the hypotheses of Definition \ref{classicaladhocdef}, then this agrees with the variant of the multiplicative Euler characteristic defined there.
\begin{proposition}\label{AgreementProp}
Consider any $A\in\cat{\zz}$, and denote by $e$ the operator $e:A\to \Sigma A$ on the underlying bounded chain complex of $\zz$-modules (which we abusively also denote by $A$).  Suppose that the cohomology groups of the chain complex $(H^*(A),e)$ are all finite.  Then the multiplicative Euler characteristic $\chi(A,e)$ from definition \ref{classicaladhocdef} agrees with the value $\chi([A],e)$ of the $K$-theory class represented by $A$ under the map in definition \ref{IntegerEulerCharGen}. 
\end{proposition}
\begin{proof}
First suppose that $A$ has finite homotopy groups, such that $A\in\cat{\zz}_{\mathrm{tors}}$.  In this case, we can use that, under the forgetful functor $\cat{\zz}\to\mo{Perf}(\zz)$, $A$ lands in the full subcategory $\mo{Perf}(\zz)_{\mathrm{tors}}$, and the induced map $K_0(\cat{\zz}_{\mathrm{tors}})\to K_0(\mo{Perf}(\zz)_{\mathrm{tors}})\to \bb{Q}^{\times}$ agrees with $K_0(\cat{\zz}_{\mathrm{tors}})\to K_0(\cat{\zz})\to\bb{Q}^{\times}$, so that $\chi(A,e)$ agrees with $\chi(A)$ by Lemma \ref{classicalagreementlemma}.

Now, suppose that $A$ is not necessarily torsion, but such that $(H^*(A),e)$ has finite cohomology groups.  Then by \cite[Proposition 3.21]{mathew2016residuefieldsclassrational}, $A\otimes_{\zz}\bb{Q}$ must decompose as a sum of shifts of $N_1$ (using that $N_i$ for $i>1$ come from an action which is inherently not semi-simple).  Picking a perfect $C(S^1,\zz)$-module $V$ and a map $V\to A$ becoming an equivalence rationally, we find that $\chi(A,e)=\chi(\mo{cofib}(V\to A),e)$, and one can check by hand that $\chi(\mo{cofib}(V\to A),e)$ agrees with the multiplicative Euler characteristic of the complex $(H^*(A),e)$.
\end{proof}
We will henceforth simply write $\chi(A,e)$ to mean $\chi([A],e)$, that is, our map applied to the $K$-theory class represented by $A$.
\subsection{Adaptation to $\ell$-adic cohomology}  We explain now how to adapt the constructions in \textsection 2.2 to $\ell$-adic \'{e}tale cohomology and syntomic cohomology, in order to define $\chi(R\Gamma(X,\zz_{\ell}(n)),e)$ in the appropriate settings.  Here we work with the $\ell$-completed variant of $\cat{\zz}$, which we denote by $\cat{\zz_{\ell}}$, defined as the full stable subcategory of $\mo{Mod}(C(S^1,\zz_{\ell}))$ generated by perfect complexes and $\fff_{\ell}$.  By Theorem \ref{KZeroComp}, we have that $$K_0(\cat{\zz_{\ell}})\simeq K_0(\cat{\zz_{\ell}}_{\tors})\oplus \zz\cdot [C(S^1,\zz_{\ell})].$$ This allows us to again define a multiplicative Euler characteristic via the group homomorphism $$\chi(-,e):K_0(\cat{\zz_{\ell}})\to K_0(\cat{\zz_{\ell}}_{\tors})\to K_0(\mo{Perf}(\zz_{\ell})_{\tors})\xrightarrow{\chi^{\times}}\bb{Q}^{\times}.$$

We consider the category $$\lge\simeq \mo{Fun}(B\zz,\mo{Perf}(\zz_{\ell}))$$ of $\ell$-adic sheaves on the \'{e}tale site of $\fff_p$ which are perfect as underlying $\zz_{\ell}$-modules, that is, the bounded derived category of lisse $\zz_{\ell}$ sheaves on a point. The main result of this subsection is the following.
\begin{theorem}\label{EtaleLifting}
There is a natural lift of the functor $$(-)^{h\zz}:\lge \to \mo{Perf}(\zz_{\ell})$$ to a functor $$(-)^{h\zz}:\lge\to \cat{\zz_{\ell}}.$$
\end{theorem}
\begin{proof}
The functor $(-)^{h\zz}$ is lax symmetric monoidal by \cite[Corollary 7.3.2.7]{HA} as it is right adjoint to the symmetric monoidal functor $$\mo{Perf}(\fff_p)\xrightarrow{(-)^{triv}}\mo{Perf}((\fff_p)_{\acute{e}t},\zz_\ell).$$ The module structure over the unit, $\zz_{\ell}$ with the trivial action, combined with this lax symmetric monoidality, gives us a canonical lift $$(-)^{h\zz}\colon\mo{Perf}((\fff_p)_{\acute{e}t},\zz_\ell)\to \mo{Mod}(C(S^1,\zz_{\ell})).$$   It suffices to show that the image lands in $\cat{\zz_{\ell}}$.

Let $A\in \mo{Perf}((\fff_p)_{\acute{e}t},\zz_\ell)$ be arbitrary.  Since specifying $A$ is the same as specifying a perfect $\zz_{\ell}$-module with a $\zz$-action, the long exact sequence attached to the fiber sequence $$A^{h\zz}\to A\xrightarrow{1-\sigma} A,$$ where $\sigma\in\zz$ is a generator, tells us that $A^{h\zz}$ is such that $\pi_*(A^{h\zz})$ is finitely generated.  By Theorem \ref{DetectionTheorem}, it suffices now to show that $A^{h\zz}\otimes_{\zz_{\ell}}\bb{Q}_{\ell}$ is perfect over $C(S^1,\bb{Q}_{\ell})$.  Note that $$A^{h\zz}\otimes_{\zz_{\ell}}\bb{Q}_{\ell}\simeq \left(A\otimes_{\zz_{\ell}}\bb{Q}_{\ell}\right)^{h\zz},$$ so we reduce to studying the functor $$(-)^{h\zz}:\mo{Perf}((\fff_p)_{\acute{e}t},\bb{Q}_\ell)\to \mo{Mod}(C(S^1,\bb{Q}_{\ell})).$$  But by Jordan decomposition, we can decompose $A\otimes_{\zz_{\ell}}\bb{Q}_{\ell}$ as $B\oplus C$, where a generator $\sigma$ of $\zz$ acts on $B$ as a generalized $\sigma=1$-eigenspace, and $1-\sigma$ is invertible on $C$.  In this case, $$(A\otimes_{\zz_{\ell}}\bb{Q}_{\ell})^{h\zz}\simeq B^{h\zz} \oplus C^{h\zz},$$ where $C^{h\zz}\simeq 0$, and $B^{h\zz}$ decomposes as a sum of shifts of modules of the form $N_m\otimes_{\zz_{\ell}}\bb{Q}_{\ell}$, which are in particular perfect, giving the claim.
\end{proof}
This leads to the following corollary, which will be crucial when proving special value formulae in the absence of Tate's semi-simplicity conjecture.
\begin{corollary}\label{CorEtaleLift}
If $X$ is a smooth proper $\fff_p$-scheme, then for any $n\in\zz$, there is a natural lift of $R\Gamma_{\acute{e}t}(X,\zz_{\ell}(n))$ to $\cat{\zz_{\ell}}$.
\end{corollary}
\begin{proof}
This follows by the previous theorem using the facts that $$R\Gamma_{\acute{e}t}(X,\zz_{\ell}(n))\simeq \left(R\Gamma_{\acute{e}t}(\overline{X},\zz_{\ell}(n))\right)^{h\zz},$$ and $$R\Gamma_{\acute{e}t}(\overline{X},\zz_{\ell}(n))\in \mo{Perf}((\fff_p)_{\acute{e}t},\zz_\ell),$$ whenever $X$ is a smooth proper $\fff_p$-scheme, e.g. by \cite[Theorem 19.1]{milne2012lectures}.
\end{proof}

\subsection{Adaptation to $p$-adic cohomology}  When we move to the world where $\ell=p$, it turns our the right system of coefficients to look at is that is no longer $p$-adic \'{e}tale cohomology, but rather that of prismatic F-gauges, arising from work of Drinfeld and Bhatt-Lurie.  We refer to \cite{bhatt2022prismatic} and \cite{mondal2025zetafunctionfgaugesspecial} for a more formal discussion of prismatic F-gauges, but for the purposes of this paper, we content ourselves with the following.

Roughly speaking, the category $\fge:=\fge(\Spec(\fff_p))$ of \textit{perfect prismatic F-gauges} (on a point) is the category of graded perfect $\zz_p$-modules $M^{\bdot}$, equipped with maps $t:M^{\bdot}\to M^{\bdot-1}$ and $u:M^{\bdot}\to M^{\bdot+1}$ such that $ut=tu=p$, equipped with equivalence $M[1/u]\simto M[1/t]$.  To be more formal about this, consider the graded ring $\zz_p[u,t]/(tu=p)\in\mo{Gr}(\mc{D}(\zz_p))$ with $u$ in degree 1 and $t$ in degree -1.  We can talk about graded perfect modules over this ring, defining a category $\mo{Perf}^{gr}(\zz_p[u,t]/(tu=p))$, which is the category of dualizable objects inside the symmetric monoidal category of graded $\zz_p[u,t]/(tu=p)$-modules.  The following makes explicit the category of perfect prismatic $F$-gauges on a point, implicitly using the formal GAGA from \cite[Lemma 3.4.11]{bhatt2022prismatic}.
\begin{definition}\label{EvilFGaugeDef}
The category $\fge$ of \textit{perfect prismatic $F$-gauges} is defined as the pullback of stable $\infty$-categories
\begin{center}
	\begin{tikzcd}
		\fge\rar\dar & \mo{Perf}^{gr}(\zz_p[u,t]/(tu=p))\dar{(-)[1/u]\oplus (-)[1/t]}\\
		\mo{Perf}(\zz_p)\rar{\text{diag}} & \mo{Perf}(\zz_p)\oplus \mo{Perf}(\zz_p).
	\end{tikzcd}
\end{center}
That is, a prismatic $F$-gauge consists of the data of a perfect graded $\zz_p[u,t]/(tu=p)$-module $M$ together with an equivalence $$\varphi:\varinjlim_{u}M^{\bdot}=M^{\infty}\simto M^{-\infty}=\varinjlim_{t}M^{\bdot},$$ which we suggestively write as a Frobenius, in $\mo{Perf}(\zz_p)$.
\end{definition}
\begin{remark}\label{ExplainEvil}
As indicated prior to the definition, this category we define is really that of perfect prismatic $F$-gauges on $\mo{Spec}(\fff_p)$, which is to say, perfect complexes on the syntomification $\fff_p^{\mo{syn}}$, with this notation taken from \cite{bhatt2022prismatic}.  For the scope of this paper, we will assume we are always working in this category, largely ignoring how we got here, so we have need only to examine these coefficients on a point.
\end{remark}
Before we prove the analogue of Theorem \ref{EtaleLifting}, let's record some basic statements about prismatic F-gauges.
\begin{definition}[\cite{bhatt2022prismatic}, Definition 4.4.1]\label{SyntomicDef}
If $N\in\fge$, we define the \textit{syntomic complex} of $N$ as $$R\Gamma_{\syn}(N)=\mo{Hom}_{\fge}(\mathbbl{1},N)\in\mc{D}(\zz_p),$$ where $\mathbbl{1}$ is the unit object in $\fge$.
\end{definition}
\begin{lemma}[\cite{bhatt2022prismatic}, Remark 4.2.8]\label{synDeftu}
Under our description from \ref{EvilFGaugeDef}, the syntomic cohomology of a given $M\in\fge$ is computed by $$R\Gamma_{\syn}(M)=\mo{fib}(M^0\xrightarrow{t^{\infty}-\varphi\circ u^{\infty}}M[1/t]),$$ where $t^{\infty}:M^0\to M^{-\infty}=M[1/t]$ and $u^{\infty}:M^0\to M^{\infty}=M[1/u]$ are the canonical inclusions.
\end{lemma}
\begin{proof}
This follows from the pullback square in the definition.
\end{proof}
\begin{lemma}[\cite{bhatt2022prismatic}, Proposition 4.5.1]\label{LemSyntomicBounded}
If $N\in\fge$, then $R\Gamma_{\syn}(N)\in \mc{D}^b(\zz_p)$.
\end{lemma}
\begin{lemma}[\cite{mondal2025zetafunctionfgaugesspecial}, Proposition 2.9]\label{LemRationalFGauge}
The rationalization $\fge[1/p]$ is equivalent to the category $\mo{Fun}(B\zz,\mo{Perf}(\bb{Q}_p))$ of group representations of $\zz$ on perfect $\bb{Q}_p$-complexes.  In particular, $R\Gamma_{\syn}(-)[1/p]$ identifies with $$(-)^{h\zz}:\mo{Fun}(B\zz,\mc{D}^b(\bb{Q}_p))\to \mc{D}^b(\bb{Q}_p).$$
\end{lemma}
\begin{lemma}[\cite{bhatt2022prismatic}, Remark 4.2.7]\label{pEtaleComparison}
There is a symmetric monoidal functor $\mo{Fun}(B\zz,\mc{D}^b(\zz_p))\to \fge$ which induces an equivalence on endomorphism objects for the unit.
\end{lemma}
\begin{proof}
This follows from the Riemann-Hilbert correspondence in \cite[Remark 4.2.7]{bhatt2022prismatic}, to which we refer for a more formal presentation, where we identify $\mo{Fun}(B\zz,\mc{D}^b(\zz_p))$ with the bounded derived category of lisse $\zz_p$-sheaves on $\fff_p$.  To give an ad libbed description of the functor, it can (roughly) be described explicitly by fixing a generator $\sigma$ of $\zz$, then proceed by taking a perfect $\zz_p$-module $M$ to $M^i=M$, $t:M^{i}\to M^{i-1}$ is the identity if $i\leq 0$, and multiplication by $p$ for $i>0$, while $u:M^{i}\to M^{i+1}$ is the identity for $i\geq 0$, and multiplication by $p$ for $i<0$, where finally we let $\varphi:M[1/u]\simeq M\simto M\simeq M[1/t]$ be the map $\sigma$.
\end{proof}
With these preliminaries in hand, we can now prove.
\begin{theorem}\label{FGaugeLift}
The functor $$R\Gamma_{\syn}(-):\fge\to \mc{D}^b(\zz_p)$$ lifts canonically to a functor $$R\Gamma_{\syn}(-):\fge\to \cat{\zz_p}.$$
\end{theorem}
\begin{proof}
Using that $R\Gamma_{\syn}$ is right adjoint to the symmetric monoidal functor $$-\otimes_{\zz_p}\zz_{p}[u,t]/(tu=p):\mc{D}^b(\zz)\to \fge,$$ it is again lax symmetric monoidal, so it lifts to a functor $$R\Gamma_{\syn}:\fge\to\mo{Mod}(R\Gamma_{\syn}(\mathbbl{1})).$$ Using Lemma \ref{pEtaleComparison}, we find that $R\Gamma_{\syn}(\mathbbl{1})\simeq C(S^1,\zz_p)$, and it remains to show that the image of $R\Gamma_{\syn}(-)$ lands in $\cat{\zz_p}$.  This in turn follows by Theorem \ref{DetectionTheorem}, using Lemma \ref{LemSyntomicBounded}, and \ref{LemRationalFGauge}, along with the fact that $$(-)^{h\zz}:\mo{Fun}(B\zz,\mc{D}^b(\bb{Q}_p))\to \mo{Mod}(C(S^1,\bb{Q}_p))$$ lands in the subcategory of perfect modules.
\end{proof}
\begin{corollary}\label{SchemeFGaugeLift}
For any smooth proper scheme $X/\Spec(\fff_p)$, there is a canonical lift of $R\Gamma_{\syn}(X,\zz_p(n))$ to an object in $\cat{\zz_p}$.
\end{corollary}
\begin{proof}
This follows by \cite[Remark 4.2.3]{bhatt2022prismatic}, which gives us that the prismatic $F$-gauge $\mc{H}_{\syn}(X)$ attached to a smooth proper $\fff_p$-scheme $X$ lands in $\fge$, combined with \ref{FGaugeLift}.
\end{proof}
\newpage
\section{Special Values of $\zeta$ Functions}
\indent \indent In this section, we will provide a linear algebraic proof of Milne's Theorem \ref{Milne1}.  To motivate our constructions, we recall how Bayer and Neukirch produce a formula for $|\zeta(X,n)|_{\ell}$ for smooth proper $\fff_p$-schemes $X$ when $\ell\neq p$, and where $\zeta(X,s)$ has no pole at $n$ (which holds e.g. for $n<0$ or $n>\mo{dim}(X)$), but we will proceed with a version of the proof due to Schneider.  We will then proceed by a similar proof to recover Schneider's results at $0\leq n\leq \mo{dim}(X)$ without the need for Tate's semi-simplicity conjecture, via a method that works for non-semi-simple Galois actions (in particular, with zeta functions appropriately defined, extending these formulae to coefficients in $\lge$, where there are many examples of Frobenius actions which fail to be semi-simple).  Finally, we will explain how to use similar linear algebraic observations to adapt this style of proof to the $\ell= p$ case.
\subsection{Schneider's proof in the $\ell\neq p$ case, with $n<0$ or $n>\mo{dim}(X)$}  Let $X$ be a smooth proper $\fff_p$-scheme, and fix a prime $\ell\neq p$.  Suppose that $n<0$ or $n>\mo{dim}(X)$, so that $\zeta(X,s)$ does not have a pole at $n$.  Recall by the proof of the Weil conjectures, \cite{grothendieck1964formule} we have that $$\zeta(X,s)=\prod_{i=0}^{2\mo{dim}(X)}\mo{det}(1-\varphi p^{-s}|H^{i}_{\acute{e}t}(\overline{X},\bb{Q}_{\ell}))^{(-1)^{i+1}},$$ where $\varphi$ denotes the Frobenius operator acting on the absolute $\ell$-adic \'{e}tale cohomology of $X$.  In order to determine the $\ell$-adic absolute value $|\zeta(X,n)|_{\ell}$, it suffices to determine $$|\mo{det}(1-\varphi p^{-n}|H^i_{\acute{e}t}(\overline{X},\bb{Q}_{\ell}))|_{\ell}$$ for each $0\leq i\leq 2\mo{dim}(X)$.  Since $\zeta(X,s)$ does not have a zero or a pole at $n$ (by the Riemann Hypothesis proven by Deligne \cite{deligne1974conjecture}), $1-\varphi p^{-n}$ is an invertible operator on each finite-dimensional $\bb{Q}_{\ell}$-vector space $H^i_{\acute{e}t}(\overline{X},\bb{Q}_{\ell})$.  The key linear algebra observation is the following
\begin{lemma}\label{LinAlgLem1}
Let $V$ be a finite dimensional $\bb{Q}_{\ell}$-vector space equipped with an automorphism $\rho:V\simto V$, and assume that there is a $\zz_{\ell}$-lattice $L\subseteq V$ with $\rho(L)\subseteq L$.  Then $$|\mo{det}(\rho\vert_{V})|_{\ell}=|L/\rho(L)|^{-1},$$ where $|L/\rho(L)|$ is the order of the finite group $L/\rho(L)$.
\end{lemma}
In our setup, we nearly have a perfect choice of distinguished $\zz_{\ell}$-lattice: integral \'{e}tale cohomology, which is finitely generated by \cite[Theorem 19.1]{milne2012lectures}.  The reason we have to say nearly is due to the possible existence of torsion.  Nevertheless, we can proceed roughly as if this were a non-issue.  The following theorem was first proved by Bayer and Neukirch in \cite{Neukirch1978/79}, but the proof we provide is adapted from Schneider \cite[Theorem 5]{Schneider1982}.
\begin{theorem}[\cite{Neukirch1978/79}, \cite{Schneider1982}]\label{Valueell1}
Let $X$ be a smooth proper $\fff_p$-scheme, $\ell\neq p$ a prime, and $n\in\zz$ such that either $n<0$ or $n>\mo{dim}(X)$.  Then, we have that
$$|\zeta(X,n)|_{\ell}^{-1}=\chi(R\Gamma_{\acute{e}t}(X,\zz_{\ell}(n))),$$ where $\zz_{\ell}(n)$ is the $n$th Tate twist of $\zz_{\ell}$.
\end{theorem}
\begin{proof}
Using equation (\ref{Zetaboi}), we have that $$|\zeta(X,n)|_{\ell}^{-1}=\prod_{i=0}^{2\mo{dim}(X)}|\mo{det}(1-\varphi p^{-n}|H^{i}_{\acute{e}t}(\overline{X},\bb{Q}_{\ell}))|_{\ell}^{(-1)^{i}}.$$  Now, by Lemma \ref{LinAlgLem1}, using the $\zz_{\ell}$-lattice $$H^{i}_{\acute{e}t}(\overline{X},\zz_{\ell})/\mathrm{tors}\subseteq H^i_{\acute{e}t}(\overline{X},\bb{Q}_{\ell}),$$ we find that 
\begin{align*}
|\zeta(X,n)|_{\ell}^{-1}&=\prod_{i=0}^{2\mo{dim}(X)}|\mo{det}(1-\varphi p^{-n}|H^{i}_{\acute{e}t}(\overline{X},\bb{Q}_{\ell}))|_{\ell}^{(-1)^{i}}\\
&=\prod_{i=0}^{2\mo{dim}(X)}\left|\left(H^{i}_{\acute{e}t}(\overline{X},\zz_{\ell})/\mathrm{tors}\right)/\left(1-\varphi p^{-n}\left(H^{i}_{\acute{e}t}(\overline{X},\zz_{\ell})/\mathrm{tors}\right)\right)\right|^{(-1)^{i+1}}.
\end{align*}
Thus, our claim reduces to showing that the RHS agrees with the multiplicative Euler characteristic $\chi(R\Gamma_{\acute{e}t}(X,\zz_{\ell}(n)))$.  Indeed, the complex $R\Gamma_{\acute{e}t}(X,\zz_{\ell}(n))$ is obtained by the homotopy fixed points of the Frobenius action on $R\Gamma_{\acute{e}t}(\overline{X},\zz_{\ell}(n)),$ where $\varphi$ acts as $\varphi p^{-n}$ would on $R\Gamma_{\acute{e}t}(\overline{X},\zz_{\ell})$.  That is to say, we have a fiber sequence 
$$ R\Gamma_{\acute{e}t}(X,\zz_{\ell}(n))\to R\Gamma_{\acute{e}t}(\overline{X},\zz_{\ell})\xrightarrow{1-\varphi p^{-n}} R\Gamma_{\acute{e}t}(\overline{X},\zz_{\ell}).$$  Looking at the induced long exact sequence on homotopy groups, we decompose $$H^i_{\acute{e}t}(\overline{X},\zz_{\ell})\simeq H^i_{\acute{e}t}(\overline{X},\zz_{\ell})/\mathrm{tors} \oplus \mathrm{tors}$$ in such a way that the torsion-free summand maps to itself under $1-\varphi p^{-n}$.  Then, $$|\ker(1-\varphi p^{-n}|\mathrm{tors})|=|\mo{coker}(1-\varphi p^{-n}|\mathrm{tors})|,$$ which implies that the torsion summand contributes groups of equal size to $H^{i}(R\Gamma_{\acute{e}t}(X,\zz_{\ell}(n)))$ and $H^{i+1}(R\Gamma_{\acute{e}t}(X,\zz_{\ell}(n)))$, canceling out the contribution to the multiplicative Euler characteristic.  Thus, the only contribution we see to $\chi(R\Gamma_{\acute{e}t}(X,\zz_{\ell}(n)))$ is from $$\mo{coker}(H^i_{\acute{e}t}(\overline{X},\zz_{\ell})/\mathrm{tors}\xrightarrow{1-\varphi p^{-n}}H^i_{\acute{e}t}(\overline{X},\zz_{\ell})/\mathrm{tors}),$$ injecting into $H^{i+1}(R\Gamma_{\acute{e}t}(X,\zz_{\ell}(n)))$, thus contributing a factor of $$\left|\left(H^{i}_{\acute{e}t}(\overline{X},\zz_{\ell})/\mathrm{tors}\right)/\left(1-\varphi p^{-n}\left(H^{i}_{\acute{e}t}(\overline{X},\zz_{\ell})/\mathrm{tors}\right)\right)\right|^{(-1)^{i+1}}$$ to $\chi(R\Gamma_{\acute{e}t}(X,\zz_{\ell}(n)))$.  Since this is the case for all $i$, we conclude that $$\chi(R\Gamma_{\acute{e}t}(X,\zz_{\ell}(n)))=\prod_{i=0}^{2\mo{dim}(X)}\left|\left(H^{i}_{\acute{e}t}(\overline{X},\zz_{\ell})/\mathrm{tors}\right)/\left(1-\varphi p^{-n}\left(H^{i}_{\acute{e}t}(\overline{X},\zz_{\ell})/\mathrm{tors}\right)\right)\right|^{(-1)^{i+1}},$$ as claimed.
\end{proof}
\subsection{The case $\ell\neq p$ with $0\leq n\leq \mo{dim}(X)$}  If $0\leq n\leq \mo{dim}(X)$, there might be a pole of $\zeta(X,s)$ at the integer $n$, so we cannot ask for a value of $\zeta(X,n)$.  Instead, we note that $\zeta(X,s)\sim C(X,n)\cdot (1-p^{n-s})^{-\rho_n}$ as $s\to n$, for $\rho_n$ the order of the pole at $n$, and $C(X,n)$ a fixed rational number, the value of which we can study.  Now, in \cite{Schneider1982}, assuming Tate's semi-simplicity conjecture, Schneider was able to show (albeit in a slightly different language)
\begin{theorem}[\cite{Schneider1982}, Theorem 5]\label{WithTateAssumptionell}
Let $X$ be a smooth proper scheme over $\fff_p$, and suppose that Tate's semi-simplicity conjecture holds for $X$.  Then, taking $e$ to be cup product with the Euler class, we have that
\begin{equation}
|C(X,n)|_{\ell}^{-1}=\chi(R\Gamma_{\acute{e}t}(X,\zz_{\ell}(n)),e).
\end{equation}
\end{theorem}
Using our extension of the multiplicative Euler characteristic to $\cat{\zz}_{\ell}$, and Theorem \ref{EtaleLifting}, we can remove the assumption of Tate's semi-simplicity conjecture, allowing us to prove.
\begin{theorem}\label{MainTheoremEll}
Let $X$ be a smooth proper scheme over $\fff_p$.  Then, considering $R\Gamma_{\acute{e}t}(X,\zz_{\ell}(n))\in\cat{\zz_{\ell}}$, we have that
\begin{equation}\label{Valueell}
	|C(X,n)|_{\ell}^{-1}=\chi(R\Gamma_{\acute{e}t}(X,\zz_{\ell}(n)),e).
\end{equation}
\end{theorem}
\begin{proof}
Our goal is to attempt to argue as in Theorem \ref{Valueell1}, except now we run into a problem: there might exist a pole of $\zeta(X,s)$ at $s=n$.  By the Riemann hypotheses proved in \cite{deligne1974conjecture}, the existence of such a pole is equivalent to the operator $1-\varphi/p^n$ failing to be invertible on the $\bb{Q}_{\ell}$-vector space $H^{2n}_{\acute{e}t}(\overline{X},\bb{Q})$.  We can decompose the generalized $\varphi=1$-eigenspace of the action of $\varphi$ on $R\Gamma_{\acute{e}t}(\overline{X},\bb{Q}_{\ell}(n))$ into Jordan blocks, which allows us to take indecomposable summands $W_{1},\ldots,W_k$ of this generalized eigenspace (which are really $\bb{Q}$-vector spaces shifted into the appropriate homological degree), which include $\zz$-equivariantly into $R\Gamma_{\acute{e}t}(\overline{X},\bb{Q}_{\ell}(n))$.  

We may choose $V_i$ a $\zz_{\ell}$-linear Jordan block for $\varphi=1$ of the same rank as $W_i$, and consider the $\zz$-equivariant map $$V_i\to W_i\to R\Gamma_{\acute{e}t}(\overline{X},\bb{Q}_{\ell}(n)).$$  Since $V_i$ is a compact object in $\lge$, and $$R\Gamma_{\acute{e}t}(\overline{X},\bb{Q}_{\ell}(n))=\colim_{ \ell}R\Gamma_{\acute{e}t}(\overline{X},\zz_{\ell}(n)),$$ we may lift this map to a map $V_i\to R\Gamma_{\acute{e}t}(\overline{X},\zz_{\ell}(n))$ which is equivalent to the inclusion of the Jordan block $W_i$ after rationalization.  Doing this for $i=1,...,k$, we get a $\zz$-equivariant map $$s\colon\bigoplus_{i=1}^{k}V_i\to R\Gamma_{\acute{e}t}(\overline{X},\zz_{\ell}(n)).$$  Let 
\begin{equation}\label{MXDefSeq}
M(\overline{X}):=\mo{cofib}(s\colon\bigoplus_{i=1}^{k}V_i\to R\Gamma_{\acute{e}t}(\overline{X},\zz_{\ell}(n)))
\end{equation} denote the cofiber of this map.  We note that $$M(\overline{X})\otimes_{\zz_{\ell}}\bb{Q}_{\ell}\simeq R\Gamma_{\acute{e}t}(\overline{X},\bb{Q}_{\ell}(n))/\bigoplus_{i=1}^{k}W_i,$$ is such that the action of Frobenius has no $\varphi=1$ fixed points, so that $1-\varphi\vert_{M(\overline{X})[1/\ell]}$ is invertible, and the same argument in Theorem \ref{Valueell1} tells us that 
$$\prod_{i\in\zz}|\mo{det}(1-\varphi\vert|_{M(\overline{X})[1/\ell]})|_{\ell}^{(-1)^{i}}=\chi(M(\overline{X})^{h\zz}).$$
Since the left hand side is equivalent to the inverse of the $\ell$-adic absolute value of $C(X,n)$, it suffices to show that $\chi(R\Gamma_{\acute{e}t}(X,\zz_{\ell}(n)),e)=\chi(M(\overline{X})^{h\zz},e)$.  Applying homotopy fixed points for $\zz$ to the defining fiber sequence (\ref{MXDefSeq}), we get a fiber sequence in $\cat{\zz}_{\ell}$:
\begin{equation}
	\bigoplus_{i=1}^{k}V_i^{h\zz}\to R\Gamma_{\acute{e}t}(X,\zz_{\ell}(n))\to M(\overline{X})^{h\zz}.
\end{equation}
Since $V_i^{h\zz}$ is a perfect $C(S^1,\zz_{\ell})$-module, $\chi(V_i^{h\zz},e)=1$, and since our multiplicative Euler characteristic is multiplicative along fiber sequences, we find that
$$\chi(R\Gamma_{\acute{e}t}(X,\zz_{\ell}(n)),e)=\chi(M(\overline{X})^{h\zz},e),$$ as desired.
\end{proof}
\begin{remark}\label{remellcoeff}
	In the previous theorem, nowhere did we explicitly use that $R\Gamma_{\acute{e}t}(X,\zz_{\ell}(n))\in\lge$ came from a smooth proper scheme $X$, and indeed this theorem works more generally (when the zeta function is defined by the formula (\ref{Zetaboi}) on rationalized cohomology groups) for any $M\in\lge$, which in particular can be chosen to have Galois action which is not semi-simple on the generalized $\varphi=p^n$-eigenspace of the rationalization.
\end{remark}
\subsection{The case $\ell=p$ with $n<0$ or $n>\mo{dim}(X)$}  We now turn our attention to mimicking Schneider's argument in the setup where $\ell=p$.  Immediately we encounter problems, since $p$-adic \'{e}tale cohomology is not well-behaved, nor does it satisfy the analogue of (\ref{Zetaboi}).  The correct version of $p$-adic cohomology to look at turns out to be crystalline cohomology, after a proof of \cite{lubkin1968p}.  This takes us to our next issue, since for $n\gg 0$, we should not expect there to exist a $\zz_p$-lattice $L\subseteq H^i_{crys}(X,\bb{Q}_p)$ which is mapped to itself under the action of $1-\varphi/p^n$.  Instead, we need to make the following modification.
\begin{lemma}\label{LinAlgLem2}
Let $V$ be a finite dimensional $\bb{Q}_p$-vector space, equipped with an invertible linear map $\rho:V\to V$.  Suppose further that we are given two $\zz_p$-lattices $L^{\prime}\subseteq L\subseteq V$ in $V$ such that $\rho(L^{\prime})\subseteq L$.  Then we have that $$|\mo{det}(\rho)|_{p}=|L/\rho(L^{\prime})|^{-1}\cdot |L/L^{\prime}|.$$
\end{lemma}
With this observation in hand, we can start to ask what the lattices $L$ and $L^{\prime}$ could be represented by in cohomology.  There turns out to be a very natural choice for these ``lattices,'' the so-called Nygaard filtration on crystalline cohomology.  We will use the variant discussed in \cite[~\textsection 8]{bhatt2019topological} (see also \cite{Illusie1979ComplexeDD}), defined by the sub-complex $$\mc{N}^{\geq n}W\Omega_{X} = p^{n-1}VW\mc{O}_X\to p^{n-2}VW\Omega^1_{X}\ldots \to VW\Omega_{X}^{n-1}\to W\Omega_X^{n}\ldots,$$ where $V$ denotes the Verschiebung operator.  There is a divided Frobenius map $\varphi_n:=\varphi/p^n:\mc{N}^{\geq n}W\Omega_X\to W\Omega_X$, giving us a natural candidate to take the place of our lattice $L^{\prime}$.  The map denoted by ``$1$'' in the following is the map induced by the filtration $\mc{N}^{\geq n}W\Omega_X\to \mc{N}^{\geq 0}W\Omega_X\simeq W\Omega_X$.  We define 
\begin{equation} \label{syndefeq}
\zz_p(n)(X):=\mo{fib}(1-\varphi_n:\mc{N}^{\geq n}W\Omega_X\to W\Omega_{X}),
\end{equation}
the \textit{syntomic complex} of $X$.  Note that by \cite{bhatt2019topological}, this agrees with the complex $W\Omega_{X,log}^n[-n]$ defined by Milne in \cite{MilneValues}.  This filtration allows us to prove the following result, originally due to Milne, using only ``linear algebra.''
\begin{theorem}[\cite{MilneValues}, Theorem 0.1]\label{ValuepMilne}
Let $X$ be a smooth proper $\fff_p$-scheme, and $n\in\zz$ such that either $n<0$ or $n>\mo{dim}(X)$.  Then, we have that
$$|\zeta(X,n)|_{p}^{-1}=\chi(R\Gamma_{syn}(X,\zz_{p}(n)))\cdot \chi(W\Omega_X/\mc{N}^{\geq n}W\Omega_X).$$
\end{theorem}
\begin{proof}
From the fiber sequence (\ref{syndefeq}) defining $\zz_p(n)(X)$, we find that we have a long exact sequence 
$$\ldots \to H^i(X,\zz_p(n))\to H^i(X,\mc{N}^{\geq n}W\Omega_{X})\xrightarrow{1-\varphi_n}H^i(X,W\Omega_X)\to H^{i+1}(X,\zz_p(n))\to\ldots.$$
Since $X/\fff_p$ is smooth and proper, its crystalline cohomology groups $H^i(X,W\Omega_X)$ are finitely generated $\zz_p$-modules, and vanish for $i\gg 0$ and $i<0$, by \cite{berthelot2015notes} (see also e.g. \cite[Remark 60.24.12]{stacks-project}).  In particular, modulo torsion, $H^i(X,W\Omega_X)$ and $H^i(X,\mc{N}^{\geq n}W\Omega_X)$ give two lattices contained in $H^i(X,W\Omega_X\otimes_{\zz_p}\bb{Q}_p)$, such that $1-\varphi_n$ induces a map between them.  The torsion submodules $H^i(X,W\Omega_X)_{\mathrm{tors}},H^i(X,\mc{N}^{\geq n}W\Omega_X)_{\mathrm{tors}}$ contribute a factor of $$(|H^i(X,W\Omega_X)_{\mathrm{tors}}|/|H^i(X,\mc{N}^{\geq n}W\Omega_X)_{\mathrm{tors}}|)^{(-1)^{i+1}}$$ to $\chi(X,\zz_p(n))$.  From the long exact sequence
$$\ldots \to H^i(X,\mc{N}^{\geq n}W\Omega_{X})\to H^i(X,W\Omega_X)\to H^{i}(X,W\Omega_X/\mc{N}^{\geq n}W\Omega_X)\to H^{i+1}(X,\mc{N}^{\geq n}W\Omega_{X})\to\ldots,$$ we find that the torsion submodules also contribute a factor of $$(|H^i(X,W\Omega_X)_{\mathrm{tors}}|/|H^i(X,\mc{N}^{\geq n}W\Omega_X)_{\mathrm{tors}}|)^{(-1)^{i}}$$ to $\chi(W\Omega_X/\mc{N}^{\geq n}W\Omega_X)$ (where now we are taking a cofiber instead of a fiber, causing this off-set by 1 in the sign of the exponent).  These two terms will cancel out in the product $\chi(X,\zz_p(n))\chi(X,W\Omega_{X}/\mc{N}^{\geq n}W\Omega_{X})$, allowing us to focus only on the torsion-free quotients of $H^i(X,W\Omega_X)$ and $H^i(X,\mc{N}^{\geq n}W\Omega_X)$.

At this point, Lemma \ref{LinAlgLem2} together with our the same style argument from Theorem \ref{Valueell1} prove the claim.
\end{proof}

To show how this relates explicitly to Milne's formula , we use the following lemma.
\begin{lemma}\label{MilneCorrectionEquiv}
The multiplicative Euler characteristic of the crystalline cohomology of $X$ modulo the $n$th step in the Nygaard filtration is equal to Milne's correction factor $p^{\chi(X,\mc{O}_X,n)}$, with $\chi(X,\mc{O}_X,n)$ defined as in (\ref{Milnecorrection}).  That is, $$\chi(W\Omega_X/\mc{N}^{\geq n}W\Omega_X)=p^{\chi(X,\mc{O}_X,n)},$$ where $$\chi(X,\mc{O}_X,n)=\sum_{i=0}^{\mo{dim}(X)}\sum_{j=0}^{n}(-1)^{i+j}(n-j)h^i(X,\Omega^j).$$
\end{lemma}
\begin{proof}
The Nygaard filtration induces a finite filtration on $W\Omega_{X}/\mc{N}^{\geq n}W\Omega_X$, whose graded pieces can be identified via \cite[Lemma 8.2]{bhatt2019topological}, with $$\mc{N}^{\geq i}W\Omega_{X}/W\Omega^{\geq i+1}W\Omega_X\simeq \tau^{\leq i}\Omega_{X/\fff_p}.$$  Thus, by Proposition \ref{trivialclaim2}, $$\chi(W\Omega_{X}/\mc{N}^{\geq n}W\Omega_X)=\prod_{i=0}^{n-1}\chi(\mc{N}^{\geq i}W\Omega_{X}/W\Omega^{\geq i+1}W\Omega_X)=\prod_{i=0}^{n-1}\chi(\tau^{\leq i}\Omega_{X/\fff_p}).$$  Each complex $\tau^{\leq i}\Omega_{X/\fff_p}$ has a finite filtration on it, induced by the Hodge filtration, with graded pieces $\Omega_{X/\fff_p}^j[-j]$ (which we abusively identify with $R\Gamma(X,\Omega_{X/\fff_p}^j[-j])$ using our conventions), giving us $$\chi(\tau^{\leq i}\Omega_{X/\fff_p})=\prod_{j=0}^{i}\chi(\Omega_{X/\fff_p}^j[-j]).$$  Since $\chi(\Omega_{X/\fff_p}^j)=\prod_{k=0}^{\mo{dim(X)}}|H^k(X,\Omega_X^j)|^{(-1)^{k+j}}$, and $|H^k(X,\Omega_X^j)|=p^{\mo{dim}_{\fff_p}(H^k(X,\Omega_X^j))}=p^{h^{k(X,\Omega^j)}}$, we can put this all together to get 
\begin{align*}
\chi(W\Omega_{X}/\mc{N}^{\geq n}W\Omega_X)&=\prod_{i=0}^{n-1}\chi(\tau^{\leq i}\Omega_{X/\fff_p})=\prod_{i=0}^{n-1}\prod_{j=0}^{i}\chi(\Omega_{X/\fff_p}^j[-j])\\
&=\prod_{i=0}^{n-1}\prod_{j=0}^{i}\prod_{k=0}^{\mo{dim(X)}}|H^k(X,\Omega_X^j)|^{(-1)^{k+j}} =\prod_{i=0}^{n-1}\prod_{j=0}^{i}\prod_{k=0}^{\mo{dim(X)}}p^{\sum_{k=0}^{\mo{dim}(X)}(-1)^{k+j}h^k(X,\Omega^j)}\\
&=p^{\sum_{i=0}^{\mo{dim}(X)}\sum_{j=0}^{n}(n-j)(-1)^{i+j}h^i(X,\Omega^j)},
\end{align*}
as desired.
\end{proof}
\subsection{A digression on prismatic F-gauges}  As opposed to Theorem \ref{MainTheoremEll}, where we can note that our argument for the \'{e}tale cohomology of a scheme $X$ adapts equally well to a general coefficient in $\lge$ (see Remark \ref{remellcoeff}), we have to be a bit more careful when $\ell=p$, especially with encoding the data of the Nygaard filtration.  This is where the theory of prismatic F-gauges from \cite{bhatt2022prismatic} enters into the picture.  Let $M\in \fge$ as in Definition \ref{EvilFGaugeDef}.  It will be convenient to make the following definition in the special case that $1-\varphi$ is fixed point free on $M[1/p]$, or equivalently that $R\Gamma_{\syn}(M)$ has torsion homotopy groups.
\begin{definition}\label{MondalDefButWorse}
Let $M\in\fge$ be such that $1-\varphi$ is an equivalence on $M[1/p]$.  Then we define the \textit{zeta value} $$\zeta(M):=\prod_{i\in\zz}\left(\mo{det}(1-\varphi |H^i_{\syn}(M[1/p]))\right)^{(-1)^{i+1}}.$$
\end{definition}
\begin{remark}\label{RemMondalZetaDef}
The above is just the value at 0 of the zeta function attached to the F-gauge $M$ defined in \cite{mondal2025zetafunctionfgaugesspecial}, and we only define it here when there is no pole of Mondal's zeta function at zero.
\end{remark}
Before proceeding, we need to identify the analogue of Milne's correction factor for a general perfect prismatic $F$-gauge.  We record the following
\begin{proposition}\label{QuotientByNygaard}
The exact functor $$\fge\to\mc{D}^b(\zz_p), M\mapsto M[1/t]/M^0$$ has image contained in torsion modules, and in particular $\chi(M[1/t]/M^0)$ is always well-defined.
\end{proposition}
\begin{proof}
The canonical inclusion $M^0\to M[1/t]$ is a map between perfect $\zz_p$-modules which is an equivalence after inverting $p$, therefore the cofiber of this map is a torsion module in $\mc{D}^b(\zz_p)$.
\end{proof}
We can now prove a slightly more general analogue of Theorem \ref{ValuepMilne}, which is a special case of a theorem due to Mondal \cite[Theorem 1.1]{mondal2025zetafunctionfgaugesspecial}.
\begin{theorem}\label{MondalTheoremButWorse}
Let $M\in\fge$ be a perfect prismatic $F$-gauge, such that $R\Gamma_{\syn}(M)$ has finite cohomology groups, so that $\chi(R\Gamma_{\syn}(M))$ is defined.  Then $$|\zeta(M)|_{p}^{-1}=\chi(R\Gamma_{\syn}(M))\cdot\chi(M[1/t]/M^0).$$
\end{theorem}
\begin{proof}
We proceed as in the proof of Theorem \ref{ValuepMilne}.  In this case, we note that the ``lattices'' in $M[1/p]$ we want to take are represented by $M^0\subseteq M[1/t]$.  The map $$1-\varphi=t^{\infty}-\phi\circ u^{\infty}:M^0\to M[1/t]$$ is well-defined between these ``lattices,'' with fiber $R\Gamma_{\syn}(M)$.  The correction factor coming from ``$|L/L^{\prime}|$'' we get out of the cofiber of the canonical inclusion of ``lattices'' is measured by the cofiber $$M[1/t]/M^0=\mo{cofib}(t^{\infty}:M[1/t]\to M^0).$$ By the same arguments as in Theorem \ref{ValuepMilne}, we have that $$|\zeta(M)|_p^{-1}=\chi(R\Gamma_{\syn}(M))\cdot \chi(M[1/t]/M^0).$$
\end{proof}

\subsection{The case $\ell=p$ and $0\leq n\leq \mo{dim}(X)$}  We now prove the analogue of Theorem \ref{MainTheoremEll} for the case $\ell=p$, using the same style proof, while keeping track of the modifications that were required for Theorem \ref{ValuepMilne}.  If $X$ satisfies Tate's semi-simplicity conjecture, the following theorem is due to Milne \cite[Theorem 0.1]{MilneValues}.  The main result of Mondal \cite[Theorem 1.1]{mondal2025zetafunctionfgaugesspecial} provides an alternative approach to the following theorem, using the stable Bockstein characteristic in place of our generalized multiplicative Euler characteristic $\chi(-,e)$.
\begin{theorem}\label{MainTheoremValuep}
Let $X$ be a smooth proper $\fff_p$-scheme, and $0\leq n\leq \mo{dim}(X)$ an integer.  Then, we have that $$|\zeta(X,n)|_{p}^{-1}=\chi(R\Gamma_{\syn}(X,\zz_{p}(n)),e)\cdot \chi(W\Omega_X/\mc{N}^{\geq n}W\Omega_X).$$
\end{theorem}
\begin{proof}
We proceed similarly to the proof of Theorem \ref{MainTheoremEll}, taking a Jordan block decomposition $\bigoplus_{i=1}^{k}W_i\to R\Gamma_{crys}(X,\bb{Q}_p(n))$ for the fixed points of Frobenius under this map.  We can construct integral prismatic F-gauges $V_i\in \fge$ which even lie in the image of the \'{e}tale realization map from Lemma \ref{pEtaleComparison} which rationalize to $W_i$.  

Being in the full sub-category of $\fge$ generated by the unit, each of these $V_i$ is a compact object.  In particular, writing $$R\Gamma_{crys}(X,\bb{Q}_p(n))=\varinjlim\left(\ldots\to R\Gamma_{\Prism}(X,\bb{Z}_p)\{n\}\xrightarrow{p}R\Gamma_{\Prism}(X,\bb{Z}_p)\{n\}\to\ldots\right),$$ we find that the composite $\bigoplus_{i=1}^{k}V_i\to R\Gamma_{crys}(X,\bb{Q}_p(n))$ factors over some map $\bigoplus_{i=1}^{k}V_i\to R\Gamma_{\Prism}(X,\bb{Z}_p)\{n\}$.  Letting $M(X)$ once again be the cofiber in perfect prismatic F-gauges of the map we just defined,
$$M(X):=\mo{cofib}\left(\bigoplus_{i=1}^{k}V_i\to R\Gamma_{\Prism}(X,\bb{Z}_p)\{n\}\right),$$ we find that, since $R\Gamma_{\syn}:\fge\to \cat{\zz_p}$ is exact, and $\chi(R\Gamma_{\syn}(V_i),e)=1$, we have $$\chi(R\Gamma_{\syn}(M(X)),e)=\chi(R\Gamma_{\syn}(X,\zz_p(n)),e),$$ the complex $R\Gamma_{\syn}(M(X))$ is finite, and $M(X)[1/p]$ is the cofiber of the inclusion of the generalized eigenspace for $\varphi=1$ on $R\Gamma_{crys}(X,\bb{Q}_p(n))$, so the same proof as in Theorem \ref{MondalTheoremButWorse} shows that $$|C(X,n)|_{p}^{-1}=\chi(R\Gamma_{\syn}(M(X)),e)\chi(M(X)[1/t]/M(X)^0).$$  Since $(-)[1/t]/(-)^{0}:\fge\to\mc{D}^b(\zz_p)_{\mo{tors}}$ is exact, and it vanishes on each $V_i$, we have also that $$\chi(M(X)[1/t]/M(X)^0)=\chi(R\Gamma_{\Prism}(X,\zz_p)\{n\}[1/t]/R\Gamma_{\Prism}(X,\zz_p)\{n\}^0)=\chi(W\Omega_X/\mc{N}^{\geq n}W\Omega_X),$$ as desired.  Since the determinant of $1-\varphi$ away from the generalized $\varphi=1$-eigenspace of $R\Gamma_{crys}(X,\bb{Q}_p(n))$ agrees with the determinant of $1-\varphi$ on $M(X)[1/p]$, the claim is shown.
\end{proof}
\begin{remark}\label{MentionRecoverMondal}
The proof above used nothing special about the $F$-gauge attached to a smooth proper $\fff_p$-scheme.  In particular, defining the zeta function attached to a perfect prismatic $F$-gauge as in \cite[Definition 3.2]{mondal2025zetafunctionfgaugesspecial}, one can run the same proof as above to provide an alternate approach to Mondal's main theorem \cite[Theorem 1.1]{mondal2025zetafunctionfgaugesspecial}.
\end{remark}

\newpage
\section{The Finite Type Case}
\subsection{The $p\neq \ell$ case}  We may ask what can happen when we want to extend beyond the case of $X$ being smooth proper.  When we look for the analogue of Theorem \ref{Valueell1}, it turns out that compactly supported \'{e}tale cohomology is all we need to use.  To be precise, recall that
\begin{lemma}[\cite{milne2012lectures}, Theorem 29.8]\label{WeilDefFT}
	For any finite type $\fff_p$-scheme $X$, $$\zeta(X,s)=\prod_{i=0}^{2\mo{dim}(X)}\left(\mo{det}\left(1-\varphi p^{-s}|H^i_{c}(\overline{X},\bb{Q}_{\ell})\right)\right)^{(-1)^{i+1}},$$ where $H^i_c(-,\bb{Q}_{\ell})$ denotes the compactly supported $\ell$-adic  \'{e}tale cohomology of $X$. 
\end{lemma}
One particular feature of this construction is that $$R\Gamma_{\acute{e}t}(\overline{X},\zz_{\ell}(n))_c\in\lge,$$ which follows e.g. from \cite[Theorem 19.1]{milne2012lectures}, and the definition of compactly supported cohomology.  Therefore, we can apply the machinery of \textsection 2 on the nose.  We now prove.
\begin{theorem}\label{GeneralizationFTell}
Let $X$ be any finite type $\fff_p$-scheme.  Then, $R\Gamma_{\acute{e}t}(X,\zz_{\ell})_{c}$ has a canonical lift to $\cat{\zz_{\ell}}$, and we can compute $$|C(X,n)|_{\ell}^{-1}=\chi(R\Gamma(X,\zz_{\ell}(n)_c),e).$$
\end{theorem}
\begin{proof}
The first claim follows from Theorem \ref{EtaleLifting} combined with our observation that $R\Gamma_{\acute{e}t}(\overline{X},\zz_{\ell})_c\in \lge$.  The claim now follows by the same proof as in Theorem \ref{MainTheoremEll} (see also Remark \ref{remellcoeff}).
\end{proof}
\begin{remark}\label{Geisserellcompare}
The above should follow with more work assuming Tate's semi-simplicity conjecture, using the theory of alterations \cite{alterations}, and, assuming the semi-simplicity conjecture, is roughly the $\ell$-adic completion of the results of Geisser in \cite{geisser2005arithmetic}.
\end{remark}

\subsection{The $p=\ell$ case}  We can try to generalize the previous subsection to the case $\ell=p$, but this becomes significantly harder.  First of all, there's not a very good natural notion of compactly supported crystalline/syntomic cohomology, which is typically constructed out of a six functor formalism.  One can work with a variant of solid quasicoherent sheaves on various stacks presenting these cohomology theories, which allows for compactly supported cohomology to be defined, but taking for example the affine line $\bb{A}^1$, the compactly supported prismatic cohomology on $\bb{A}^1$ defined in this way is infinite dimensional, even rationally.  One fix for this is to use a rational p-adic cohomology theory called rigid cohomology (cf \cite{le2007rigid}), which provides a good theory of compactly supported cohomology, but has no good integral version satisfying \'{e}tale descent \cite{abe2022integralpadiccohomologytheories}, so is wholly insufficient for our purposes (though there are proposals for such theories which do not satisfy \'{e}tale descent) cf \cite{ertl2025integralpadiccohomologytheories} under the assumption of resolution of singularities, and \cite{merici2025motivicintegralpadiccohomology} without it.\footnote{Although it appears both theories use the cdh topology anyways.}

The solution we opt for in this paper is to define compactly supported cohomology not by way of six functor formalisms, but by explicitly enforcing the properties we want to hold.  More precisely, consider a finite type $\bb{F}_p$-scheme $U$, choose an open immersion $U\to X$ into a proper $\bb{F}_p$-scheme $X$ (that is, a compactification), with closed complement $Z$.  If our cohomology theory is given as a presheaf $\mc{F}:\mo{Sch}_{\bb{F}_p}^{\mathrm{ft},\mathrm{op}}\to \mc{C}$ on finite type $\fff_p$-schemes, valued in a stable $\infty$-category $\mc{C}$, we would like to be able to define the compactly supported cohomology of $U$ as:
\begin{equation}\label{compactcohdef}
	\mc{F}_c(U):=\mo{fib}(\mc{F}(X)\to \mc{F}(Z)).
\end{equation}
Attempting to run this as is, for instance taking $\mc{F}$ to be rational crystalline cohomology $R\Gamma_{crys}(-/\bb{Q}_p)$, this compactly supported cohomology is not well-defined, in part since $\mc{F}$ is not nil-invariant, and we never specified $Z$ had to be in the reduced subscheme structure.  There are more interesting reasons for which it can fail to be well-defined also.

In order to make a well-defined theory, we need to ensure that our compactly supported cohomology is independent of a choice of compactification.  This motivates us to introduce a condition which will ensure such independence: abstract blowup excision.  Recall the following definitions due to Voevodsky \cite{voevodsky2008homotopytheorysimplicialpresheaves}.
\begin{definition}\label{absblowup}
We say that an \textit{abstract blowup square} is a pullback square of finite type $\bb{F}_p$-schemes
\begin{equation}\label{blowupsquare}
\begin{tikzcd}
	Z^{\prime}\rar{i^{\prime}}\dar{q}& X^{\prime}\dar{p}\\
	Z\rar{i} & X
\end{tikzcd}
\end{equation}
with $p$ is proper, $i$ a closed immersion, and such that $p$ induces an isomorphism $X^{\prime}\bs Z^{\prime}\simto X\bs Z$.
\end{definition}
\begin{definition}\label{absblowupexcision}
Consider a presheaf $\mc{F}:\mo{Sch}^{\mathrm{ft},\mathrm{op}}_{\fff_p}\to \mc{C}$ landing in a stable $\infty$-category $\mc{C}$.  We say that $\mc{F}$ satisfies \textit{abstract blowup excision} if it takes any abstract blowup square (\ref{blowupsquare}) to a cartesian square in $\mc{C}$.
\end{definition}
The following proposition tells us that presheaves satisfying abstract blowup excision allow us to make our desired definition (\ref{compactcohdef}) in the desired manner.
\begin{proposition}\label{compactwelldef}
Suppose that $\mc{F}:\mo{Sch}^{\mathrm{ft},\mathrm{op}}_{\fff_p}\to \mc{C}$ is a presheaf valued in a stable $\infty$-category $\mc{C}$ which satisfies abstract blowup excision.  Take a finite type $\fff_p$-scheme $U$, together with a compactification $U\to X$ with closed complement $Z$.  Then the isomorphism class of the object $$\mc{F}_c(U):=\mo{fib}(\mc{F}(X)\to \mc{F}(Z))$$ depends only on $U$.
\end{proposition}
\begin{proof}
Suppose we have two compactifications $U\to X$ and $U\to X^{\prime}$ of our given scheme $U$.  Replacing $XT{\prime}$ by the closure of $U$ mapping in via the diagonal to $X\times X^{\prime}$ if necessary, we may assume without loss of generality that $U\to X$ factors as $U\to X^{\prime}\to X$.  Let $Z=X\bs U$, and $Z^{\prime}=Z\times_{X}X^{\prime}=X^{\prime}\bs U$.  The map $X^{\prime}\to X$ is proper, and $X^{\prime}\bs Z^{\prime}\simeq U\simeq X\bs Z$.  That is to say, we have an abstract blowup square 
\begin{center}
	\begin{tikzcd}
		Z^{\prime}\dar\rar & X^{\prime}\dar\\
		Z\rar& X.
	\end{tikzcd}
\end{center}
Since $\mc{F}$ was assumed to satisfy abstract blowup excision, applying it to this abstract blowup square yields a cartesian square in $\mc{C}$
\begin{center}
	\begin{tikzcd}
		\mc{F}(X)\rar\dar & \mc{F}(Z)\dar\\
		\mc{F}(X^{\prime})\rar & \mc{F}(Z^{\prime}).
	\end{tikzcd}
\end{center}
A cartesian square in a stable $\infty$-category induces an equivalence on the fibers of the horizontal morphisms, inducing an equivalence $\mo{fib}(\mc{F}(X)\to \mc{F}(Z))\simto \mo{fib}(\mc{F}(X^{\prime})\to\mc{F}(Z^{\prime}))$.  This shows that the compactly supported cohomology of $U$ defined with respect to $X$ is isomorphic to the compactly supported cohomology of $U$ defined with respect to $X^{\prime}$, proving the claim.
\end{proof}

This tells us what kinds of presheaves we can define compactly supported cohomology for, but we are interested in defining this for cohomology theories which expressly do not satisfy abstract blowup excision.  In order to rectify this issue, we must change the value of our cohomology theory at least on some schemes, and one way we go about this is through sheafification with respect to a topology whose definition (due to Suslin-Voevodsky \cite{suslin2000relative}) we now recall.
\begin{definition}\label{cdhtopology}
The \textit{cdh topology} on $\mo{Sch}^{\mathrm{ft}}_{\fff_p}$ is the Grothendieck topology generated by Nisnevich covers and covers of the form $\{Z\coprod X^{\prime}\to X\}$ whenever
\begin{center}
	\begin{tikzcd}
		Z^{\prime}\dar\rar & X^{\prime}\dar\\
		Z\rar& X
	\end{tikzcd}
\end{center}
is an abstract blowup square.  Given a Nisnevich sheaf $\mc{F}:\mo{Sch}^{\mathrm{ft},\mathrm{op}}_{\fff_p}\to \mc{C}$ valued in a stable $\infty$-category $\mc{C}$, we will write $L_{cdh}\mc{F}$ for the sheafification of $\mc{F}$ with respect to the cdh topology.
\end{definition}
This leads to the following observation, which in our setup is due to Voevodsky (see also \cite[Proposition 2.1.5]{elmanto2020cdhdescentcdarcdescent} for a precise reference and generalization):
\begin{proposition}[\cite{voevodsky2008homotopytheorysimplicialpresheaves}]\label{cdhsheafprop}
A Nisnevich sheaf $\mc{F}:\mo{Sch}^{\mathrm{ft},\mathrm{op}}_{\fff_p}\to \mc{C}$ valued in a stable $\infty$-category $\mc{C}$ is a cdh sheaf if and only if $\mc{F}$ satisfies abstract blowup excision in the sense of Definition \ref{absblowupexcision}.
\end{proposition}

Now, before relating this back to special values, we make a few simple observations which will be crucial in proving the main theorem of this subsection.

\begin{proposition}\label{trivialclaimzeta1}
Suppose $U$ is a finite type $\bb{F}_p$-scheme, and we are given a compactification $U\to X$ with closed complement $Z$.  Then $\zeta(X,s)=\zeta(U,s)\cdot \zeta(Z,s)$.
\end{proposition}
\begin{proof}
Recall that $\zeta(X,s)$ is the holomorphic extension of the product series $\prod_{x\in X_{(0)}}\f{1}{1-|\kappa(x)|^{-s}}$ where this converges, where $X_{(0)}$ is the set of closed points in $X$, and $|\kappa(x)|$ is the size of the residue field $\kappa(x)$ at the point $x$.  Since the open-closed decomposition $U\to X \leftarrow Z$ induces a decomposition $X_{(0)}=U_{(0)}\coprod Z_{(0)}$, the claim follows.
\end{proof}

\begin{proposition}\label{trivialclaimzeta2}
Suppose we are given an abstract blowup square
\begin{center}
	\begin{tikzcd}
		Z^{\prime}\dar\rar & X^{\prime}\dar\\
		Z\rar& X
	\end{tikzcd}
\end{center}
of finite type $\fff_p$-schemes.  Then $$\zeta(Z,s)\cdot \zeta(X^{\prime},s)=\zeta(X,s)\cdot\zeta(Z^{\prime},s).$$
\end{proposition}
\begin{proof}
By Proposition \ref{trivialclaimzeta2}, we note $$\zeta(X^{\prime}\bs Z^{\prime},s)=\zeta(X^{\prime},s)/\zeta(Z^{\prime},s),$$ and $$\zeta(X\bs Z,s)=\zeta(X,s)/\zeta(Z,s).$$ Since $X^{\prime}\bs Z^{\prime}\simeq X\bs Z$, the result follows.
\end{proof}
We record one last observation to show that the cdh sheafification of a sheaf $\mc{F}:\mo{Sch}^{\mathrm{ft},\mathrm{op}}_{\fff_p}\to \cat{R}$ is not much different than the cdh sheafification of $\mc{F}:\mo{Sch}^{\mathrm{ft},\mathrm{op}}_{\fff_p}\to \mc{D}^b(R)$ for a Dedekind ring $R$ (such as $R=\zz$, $\zz_{\ell}$).  In particular, if a map from a Nisnevich sheaf to its cdh sheafification considered as a functor to $\mc{D}^b(R)$ is an equivalence on smooth (or smooth proper) schemes, then this remains true for the cdh sheafification of our sheaf considered as landing in $\cat{R}$.  This will follow from the following general claim.
\begin{proposition}\label{basicsitethings}
	Let $\mc{F}:\mc{C}^{op}\to \cat{R}$ be a presheaf on a site $(\mc{C},\tau)$, valued in the category of Arithmetic $R$-modules for some Dedekind ring $R$, and let $i:\cat{R}\to\mc{D}^b(R)$ be the forgetful functor.  Then, letting $L_{\tau}(-)$ be the $\tau$-sheafification functor, we have $i\circ L_{\tau}\mc{F}\simeq L_{\tau}(i\circ \mc{F})$.
\end{proposition}
\begin{proof}
Since $i:\cat{R}\to \mc{D}^b(R)$ preserves all limits and colimits, and detects equivalences, it follows that $i\circ L_{\tau}(\mc{F})$ is a $\tau$-sheaf.  The right adjoint to this functor is given by taking a $\tau$-sheaf $\mc{G}:\mc{C}^{op}\to \mc{D}^b(R)$ to the presheaf $L_{\tau}\mo{Hom}_{R}(C(S^1,R),\mc{G}(-))$ valued in $\cat{R}$.  Since $C(S^1,R)\simeq R\oplus \Sigma^{-1}(R)$ as an $R$-module, we have $$i(\mo{Hom}_{R}(C(S^1,R),\mc{G}(-)))\simeq \mc{G}\oplus \Sigma \mc{G},$$ which is a finite colimit of sheaves.  Since our sheaves are valued in a stable $\infty$-category, this is already a $\tau$-sheaf, so that $$L_{\tau}\mo{Hom}_{R}(C(S^1,R),\mc{G}(-))\simeq \mo{Hom}_{R}(C(S^1,R),\mc{G}(-)),$$ and this functor is equivalent to first forgetting $\mc{G}$ to a $\mc{D}^b(R)$-valued presheaf and then applying $\mo{Hom}_{R}(\underline{C(S^1,R)},-)$, which itself is right adjoint to $L_{\tau}(i(-))$.  Since any two left adjoints to a given functor are equivalent, we get an equivalence $L_{\tau}(i(-))\simeq i(L_{\tau}(-))$.
\end{proof}

We now come to the first theorem of this section. 
\begin{theorem}\label{MainThFT}
Let $\mc{G}\subseteq \mo{Sch}^{\mathrm{prop}}_{\fff_p}$ be the smallest subcategory of proper $\fff_p$-schemes\footnote{Which we take by convention to include the empty scheme.} which contains all smooth proper $\fff_p$-schemes $X$ for which all of the maps induced by sheafification
\begin{equation}\label{cdhequivs}
\bb{Z}_p(n)(X)\to (L_{cdh}\zz_p(n))(X)\quad \text{ and }\quad \mo{R}\Gamma(X,W\Omega/\mc{N}^{\geq k}W\Omega)\to L_{cdh}\mo{R}\Gamma(X,W\Omega/\mc{N}^{\geq k}W\Omega),
\end{equation}
are equivalences, and which has the property that given an abstract blowup square 
\begin{center}
	\begin{tikzcd}
		Z^{\prime}\dar\rar & X^{\prime}\dar\\
		Z\rar& X
	\end{tikzcd}
\end{center}
with any three of the displayed schemes in $\mc{G}$, then the fourth is as well.

Then, for any finite type $\fff_p$-scheme $U$ having a compactification $U\to X$, with closed complement $Z$, such that $X,Z\in \mc{G}$, $\chi((L_{cdh}\zz_p(n))_c(U),e)$ is well-defined, and we have 
\begin{equation}\label{FTValueformula}
	|C(U,n)|_p^{-1}=\chi((L_{cdh}\zz_p(n))_c(U),e)\chi((L_{cdh}W\Omega_{-}/\mc{N}^{\geq n}W\Omega_{-})_c(U)).
\end{equation}
\end{theorem}
\begin{proof}
To begin, we prove that (\ref{FTValueformula}) holds for all $X\in\mc{G}$, where properness ensures that the compactly supported cohomology of $X$ agrees with the value on $X$ of the cdh-sheafification of our cohomology theory.  Note that for any smooth proper scheme $X$ such that the maps (\ref{cdhequivs}) are equivalences, the claim follows by Theorem \ref{MainTheoremValuep}.

Consider the class of proper $\fff_p$-schemes for which we have the formula (\ref{FTValueformula}), and such that $(L_{cdh}\zz_p(n))(X)\in\cat{\zz_p}$ lives in $\cat{\zz_p}$.  We claim that this class is closed under 3-out-of-4 for abstract blowup squares, which combined with our observation that it contains smooth proper $\fff_p$-schemes such that the maps (\ref{cdhequivs}) are equivalences, will show it contains all of $\mc{G}$.  Indeed, take an abstract blowup square 
\begin{center}
	\begin{tikzcd}
		Z^{\prime}\dar\rar & X^{\prime}\dar\\
		Z\rar& X
	\end{tikzcd}
\end{center}
with any three of the schemes satisfying the formula (\ref{FTValueformula}).  By Proposition \ref{trivialclaimzeta2}, we have that $$\zeta(Z,s)\cdot \zeta(X^{\prime},s)=\zeta(Z^{\prime},s)\cdot \zeta(X,s),$$ and then so too with $\zeta(-,s)$ replaced by $\zeta(-,n)$ or $C(-,n)$.  For any cdh sheaf $\mc{F}$, such as our given cdh sheafifications of interest, the cdh excision condition takes our abstract blowup square to a cartesian square.  This induces a fiber sequence $$\mc{F}(X)\to \mc{F}(Z)\oplus \mc{F}(X^{\prime})\to \mc{F}(Z^{\prime}).$$  If $\mc{F}:\mo{Sch}^{\mathrm{ft},\mathrm{op}}_{\fff_p}\to\mc{D}(\zz_p)$ is such that three out of the four sheaves pictured above have finite cohomology groups, so does the fourth, and $$\chi(\mc{F}(X))\cdot \chi(\mc{F}(Z^{\prime}))=\chi(\mc{F}(X^{\prime}))\cdot \chi(\mc{F}(Z)),$$ giving the correction factor part of (\ref{FTValueformula}).  If $\mc{F}:\mo{Sch}^{\mathrm{ft},\mathrm{op}}\to\mo{Mod}(C(S^1,\zz_p))$ is such that three out of the four $C(S^1,\zz_p)$-modules we get from an abstract blowup square land in $\cat{\zz_p}$, then so does the fourth, as $\cat{\zz_p}$ is stable.  In this case, the fiber sequence also gives us that $$\chi(\mc{F}(Z),e)\cdot \chi(\mc{F}(X^{\prime}),e)=\chi(\mc{F}(X),e)\cdot \chi(\mc{F}(Z^{\prime}),e).$$  We find that, as these two product formulas agree with the formula we get from Proposition \ref{trivialclaimzeta2} for the zeta functions, and the $C(X,n)$ are definitionally nonzero, the fact that three out of four schemes satisfy (\ref{FTValueformula}) implies the fourth does as well, proving the claim for every $X\in\mc{G}$.

Finally, we conclude by following essentially the exact same argument with fiber sequences.  Take for any $U$ as in the statement, any chosen compactification $U\to X$ with closed complement $Z$ such that $X,Z\in\mc{G}$.  For the two cdh sheaves $\mc{F}$ we consider, the defining fiber sequence $$\mc{F}_c(U)\to \mc{F}(X)\to \mc{F}(Z)$$ tells us that $$\chi((L_{cdh}\zz_p(n))_c(U),e)=\chi(L_{cdh}\zz_p(n)(X),e)\cdot \chi(L_{cdh}\zz_p(n)(Z),e)^{-1},$$ and $$\chi((L_{cdh}W\Omega_{-}/\mc{N}^{\geq n}W\Omega_{-})_c(U))=\chi(L_{cdh}W\Omega_{-}/\mc{N}^{\geq n}W\Omega_{-}(X))\cdot\chi(L_{cdh}W\Omega_{-}/\mc{N}^{\geq n}W\Omega_{-}(Z))^{-1}.$$ Then using the formula we derive from Proposition \ref{trivialclaimzeta1}, $$|C(U,n)|_{p}^{-1}=|C(X,n)|_{p}^{-1}|C(Z,n)|_{p},$$ combined with our previous observation that (\ref{FTValueformula}) holds for $X$ and $Z$, proves the claim.
\end{proof}
In \cite{geisser2005arithmetic}, Geisser essentially proves that the formula (\ref{FTValueformula}) holds in general under a strong resolution of singularities assumption, together with an assumption of Tate's semi-simplicity conjecture (although Geisser works globally with the Weil-\'{e}tale topology as opposed to prime by prime).  The resolution of singularities property is the following.
\begin{definition}[\cite{geisser2005arithmetic}, Definition 2.4]\label{GeisserRes}
We say that strong resolution of singularities up to dimension $d$ holds, denoted R(d), if:
\begin{enumerate}
	\item Every finite type integral separated scheme $X$ over $\fff_p$ of dimension $\leq d$ admits a proper birational map $Y\to X$ with $Y$ smooth (necessarily of the same dimension as $X$).
	\item Every proper birational map between smooth schemes of dimension $\leq d$ is refined by a sequence of blowups along smooth centers.
\end{enumerate}
\end{definition}
 We now show how to get an analogous result under similar assumptions at the prime $p$.  If the conjectures laid out in the introduction of \cite{geisser2005arithmetic} are true for schemes of dimension $\leq d$, the following recovers Geisser's result when completed at the prime $p$.
\begin{corollary}\label{GeisserTheoremCor}
Assume that R(d) holds.  Then for any finite type $\fff_p$-scheme $U$ of dimension $\mo{dim}(U)\leq d$, we have 
$$
	|C(U,n)|_p^{-1}=\chi((L_{cdh}\zz_p(n))_c(U),e)\chi((L_{cdh}W\Omega_{-}/\mc{N}^{\geq n}W\Omega_{-})_c(U)).
$$
\end{corollary}
\begin{proof}
We proceed by induction on dimension, noting that the claim holds easily for zero dimensional schemes where $L_{cdh}(\mc{F})(X)=\mc{F}(X^{red})$.  Suppose we are given any proper $\fff_p$-scheme $X$ of dimension $\leq d$.  Geisser proves \cite[Theorem 4.7]{geisser2005arithmetic} that under R(d), the comparison maps (\ref{cdhequivs}) on the de Rham version are equivalences for all smooth proper schemes $X$ of dimension $\leq d$, which passes to the de Rham-Witt case modulo the Nygaard filtration by using the finite filtration on it induced by the Nygaard filtration, and then the Hodge filtration on the filtration quotients. Since the syntomic version of the comparison map is an equivalence for all smooth proper schemes $X$ unconditionally by \cite[Corollary 6.5]{elmanto2023motivic}, $\mc{G}$ contains all smooth proper schemes of dimension $\leq d$.  We can take irreducible components $X_1,\ldots,X_n$ of $X$ and form the abstract blowup square
\begin{center}
	\begin{tikzcd}
		Z^{\prime}\rar\dar & \coprod_{i=1}^{n}X_i\dar\\
		Z\rar & X,
	\end{tikzcd}
\end{center}
where $Z=\cup_{i\neq j}X_i\cap X_j\subseteq X$ has dimension strictly smaller than $d$.  By induction, $Z,Z^{\prime}\in\mc{G}$, so $X\in\mc{G}$ if and only if $X_i\in\mc{G}$ for all $i$, and thus we may assume without loss of generality that $X$ is irreducible.  By using the cdh cover $X\to X^{\mathrm{red}}$, we may also assume that $X$ is reduced, hence integral.  Now, choose some smooth proper $\fff_p$-scheme $Y$ of dimension $\mo{dim}(X)$ which is birational to $X$, existing by our assumption that R(d) holds.  By picking some explicit opens $U\subseteq X$ isomorphic to $V\subseteq Y$ and letting $W$ be the closure of the image of the diagonal map $U\xrightarrow{\Delta}X\times Y$, the following abstract blowup square
\begin{center}
	\begin{tikzcd}
		Z\rar\dar & W\dar\\
		Z^{\prime}\rar & Y,
	\end{tikzcd}
\end{center}
shows that $W\in \mc{G}$, using our inductive hypothesis since $Z$, $Z^{\prime}$ have dimension $<\mo{dim}(X)$, and all smooth proper schemes $Y$ are in $\mc{G}$.  Now, the abstract blowup square
\begin{center}
	\begin{tikzcd}
		Z\rar\dar & W\dar\\
		Z^{\prime\prime}\rar & X,
	\end{tikzcd}
\end{center}
together with induction again, shows that $X\in\mc{G}$ as well.

Now, for any finite type $\fff_p$-scheme $U$ of dimension $\leq d$, choose a compactification $U\to X$ with $X$ of the same dimension as $U$.  Then $X\bs U$ and $X$ are both in $\mc{G}$, and the result follows by Theorem \ref{MainThFT}.
\end{proof}

\begin{remark}\label{EMWork}
As mentioned above, the work of Elmanto-Morrow \cite[Corollary 6.5]{elmanto2023motivic} shows that $\zz_p(n)(X)\simeq L_{cdh}\zz_p(n)(X)$ for any smooth (in particular smooth proper) $\fff_p$-scheme $X$ without any assumption on resolution of singularities.  This means that, up to the de Rham factor, we almost know that the class $\mc{G}$ of proper $\fff_p$-schemes where special value formulas ``work well'' contains every smooth proper scheme.
\end{remark}

While we do have strong resolution of singularities in low dimensions, the question of resolving singularities in characteristic $p$ remains an open question in general.  Nevertheless, recent advances in motivic homotopy theory provide a suitable replacement for Hodge cohomology which exists as some cdh sheaf, and agrees with Hodge cohomology on smooth \textit{projective} schemes.  Namely, we use the theory of $\bb{A}^1$-colocalization studied recently by Annala-Hoyois-Iwasa in \cite{annala2024atiyahdualitymotivicspectra}, which was adapted to the Hodge cohomology setup by Annala-Pstragowski in \cite{annala2025noteweightfiltrationscharacteristic}.  The main results we use are as follows.
\begin{lemma}\label{cdhfromA1}
An $\bb{A}^1$-invariant motivic spectrum extends to a cdh sheaf on finite type $\fff_p$-schemes.
\end{lemma}
\begin{proof}
This is due to Cisinski \cite{Cisinski_2013}, see also Khan \cite{Khan_2024}.
\end{proof}
\begin{proposition}\label{ToniPiotr}
There is a cdh sheaf $$(\Omega^j)^{\dagger}(-):\mo{Sch}^{\mathrm{ft},\mathrm{op}}_{\fff_p}\to \mc{D}(\zz_p),$$ equipped with a map $$(\Omega^j)^{\dagger}(-)\to \Omega^j(-),$$ such that for every smooth projective $\fff_p$-scheme $X$, $$(\Omega^j)^{\dagger}(X)\to \Omega^j(X)$$ is an equivalence.\footnote{Recall we are abusively identifying $\Omega^j(X)$ with $R\Gamma(X,\Omega^j)$.}
\end{proposition}
\begin{proof}
From the construction of the Hodge-filtered de Rham cohomology as a motivic spectrum by Annala-Pstragowski in \cite[Definition 3.13]{annala2025noteweightfiltrationscharacteristic}, which inherits a module structure over $\mo{kgl}$, one can apply the $\bb{A}^1$-colocalization procedure of Annala-Hoyois-Iwasa \cite{annala2024atiyahdualitymotivicspectra}, to get an $\bb{A}^1$-invariant motivic cohomology theory representing Hodge-filtered de Rham cohomology (which was studied in \cite{annala2025noteweightfiltrationscharacteristic}).  Using Lemma \ref{cdhfromA1}, and \cite[Proposition 6.21]{annala2024atiyahdualitymotivicspectra}, we see that the cdh-local cohomology theory on finite type $\fff_p$-schemes defined by the $\bb{A}^1$-invariant motivic sheaves $(\Omega^j)^{\dagger}$ agrees with the Hodge cohomology complex $R\Gamma(X,\Omega^j)$ on smooth projective $\fff_p$-schemes $X$.
\end{proof}
Combined with the result of Elmanto-Morrow, this allows us to prove the following theorem, which extends the results for special values of zeta functions to any finite type $\fff_p$-scheme which has projective compactification with projective complement (and even a little more).
\begin{theorem}\label{Strongestuncondtheorem}
Let $\mc{G}^{\prime}\subseteq \mo{Sch}^{\mathrm{prop}}_{\fff_p}$ be the smallest full subcategory of proper $\fff_p$-schemes which contains all smooth proper $\fff_p$-schemes $X$ such that $$R\Gamma(X,(\Omega^j)^{\dagger}) \to R\Gamma(X,\Omega^j)$$ is an equivalence for all $j$, and which is closed under the three-out-of-four for abstract blowup squares (\ref{absblowup}).  Then for any finite type $\fff_p$-scheme $U$ which has a compactification $U\to X$ with closed complement $Z$ such that $X,Z\in\mc{G}^{\prime}$, we have 
$$
|C(U,n)|_p^{-1}=\chi((L_{cdh}\zz_p(n))_c(U),e)\cdot p^{\chi(U,\mc{O}_U,n)_c^{\dagger}},
$$
where $$\chi(U,\mc{O}_U,n)_c^{\dagger}:=\sum_{j=0}^{n}(-1)^{j}(n-j)\chi(R\Gamma(U,(\Omega^j)^{\dagger})_c).$$
In particular, if we can find a smooth projective compactification of $U$ with smooth projective complement, this formula applies.
\end{theorem}
\begin{proof}
This follows by the exact same proof as Theorem \ref{MainThFT}, using Lemma \ref{MilneCorrectionEquiv}, noting now that given an abstract blowup square (\ref{absblowup}), $$\chi(X,\mc{O}_X,n)^{\dagger}+\chi(Z^{\prime},\mc{O}_{Z^{\prime}},n)^{\dagger}=\chi(X^{\prime},\mc{O}_{X^{\prime}},n)^{\dagger}+\chi(Z,\mc{O}_{Z},n)^{\dagger},$$ and using once again the result \cite[Corollary 6.5]{elmanto2023motivic} that $\zz_p(n)(X)\simeq L_{cdh}\zz_p(n)(X)$ for smooth proper $\fff_p$-schemes $X$ to remove the first hypothesis from Theorem \ref{MainThFT}.
\end{proof}

\newpage
\printbibliography

\end{document}